\documentclass[a4paper, 10pt, DIV=10, headinclude=false, footinclude=false]{scrartcl}

\usepackage{amsthm, amsmath, amssymb, amsfonts}
\usepackage[utf8]{inputenc}
\usepackage[english]{babel}
\usepackage{url}
\usepackage{mathtools}
\usepackage{comment}

\usepackage{bm}
\usepackage{tikz}
\usepackage{booktabs}

\usepackage{pgfplots}
\usepackage{pgfplotstable}
\usepackage{stackengine}

\numberwithin{figure}{section}
\numberwithin{table}{section}
\numberwithin{equation}{section}

\newenvironment{abstr}[1]{ \vspace{.05in}\footnotesize
	\parindent .2in
	{\upshape\bfseries #1. }\ignorespaces}{\par\vspace{.1in}}

\newenvironment{keywords}{\begin{abstr}{Key words}}{\end{abstr}}

\DeclareMathAlphabet\mathbfcal{OMS}{cmsy}{b}{n}

\newtheorem{theorem}{Theorem}[section]
\newtheorem{lemma}[theorem]{Lemma}

\newtheorem{assumption}[theorem]{Assumption}

\theoremstyle{definition}
\newtheorem{definition}[theorem]{Definition}

\newtheorem{example}[theorem]{Example}

\newcommand{\RR}{\mathbb{R}}
\newcommand{\NN}{\mathcal{N}}

\newcommand{\II}{\mathcal{I}}

\newcommand{\hatVV}{\mathbf{\hat V}}
\newcommand{\VV}{\mathbf{V}}
\newcommand{\WW}{\mathbf{W}}
\newcommand{\uu}{\mathbf{u}}
\newcommand{\vv}{\mathbf{v}}
\newcommand{\ww}{\mathbf{w}}
\newcommand{\KK}{\mathbf{K}}
\newcommand{\ff}{\mathbf{f}}

\renewcommand{\gg}{\mathbf{g}}
\newcommand{\LL}{\mathbf{L}}
\newcommand{\MM}{\mathbf{M}}
\newcommand{\PP}{\mathbf{P}}
\newcommand{\QQ}{\mathbf{Q}}
\newcommand{\Rdd}{\mathbf{R}}
\newcommand{\III}{\mathbfcal{I}}
\newcommand{\vvarphi}{\bm{\varphi}}

\newcommand{\norm}[1]{{\left\vert\kern-0.25ex\left\vert\kern-0.25ex\left\vert #1 
    \right\vert\kern-0.25ex\right\vert\kern-0.25ex\right\vert}}

\allowdisplaybreaks[4]

\begin{document}
	
	\title{Numerical homogenization of spatial network models}
	\author{F. Edelvik$^1$, M. G\"{o}rtz$^1$, F. Hellman$^2$, G. Kettil$^1$, A. M\aa lqvist$^2$}
	
	\date{}
\maketitle
\footnotetext[1]{Fraunhofer-Chalmers Centre, Chalmers Science Park, 412 88 G\"{o}teborg, Sweden}
\footnotetext[2]{Department of Mathematical Sciences, Chalmers University of Technology and University of Gothenburg, 412 96 G\"oteborg, Sweden}
\begin{abstract}
We present and analyze a methodology for numerical homogenization of spatial networks, modelling e.g.~diffusion processes and deformation of mechanical structures. The aim is to construct an accurate coarse model of the network. By solving decoupled problems on local subgraphs we construct a low dimensional subspace of the solution space with good approximation properties. The coarse model of the network is expressed by a Galerkin formulation and can be used to perform simulations with different source and boundary data at a low computational cost. We prove optimal convergence of the proposed method under mild assumptions on the homogeneity, connectivity, and locality of the network on the coarse scale. The theoretical findings are numerically confirmed for both scalar-valued (heat conduction) and vector-valued (structural) models. 
\end{abstract}

\begin{keywords}
algebraic connectivity, discrete model, multiscale method, network model, localized orthogonal decomposition, upscaling
\end{keywords}

\section{Introduction}
\begin{figure}
  \centering
  \includegraphics[width = 0.25\textwidth]{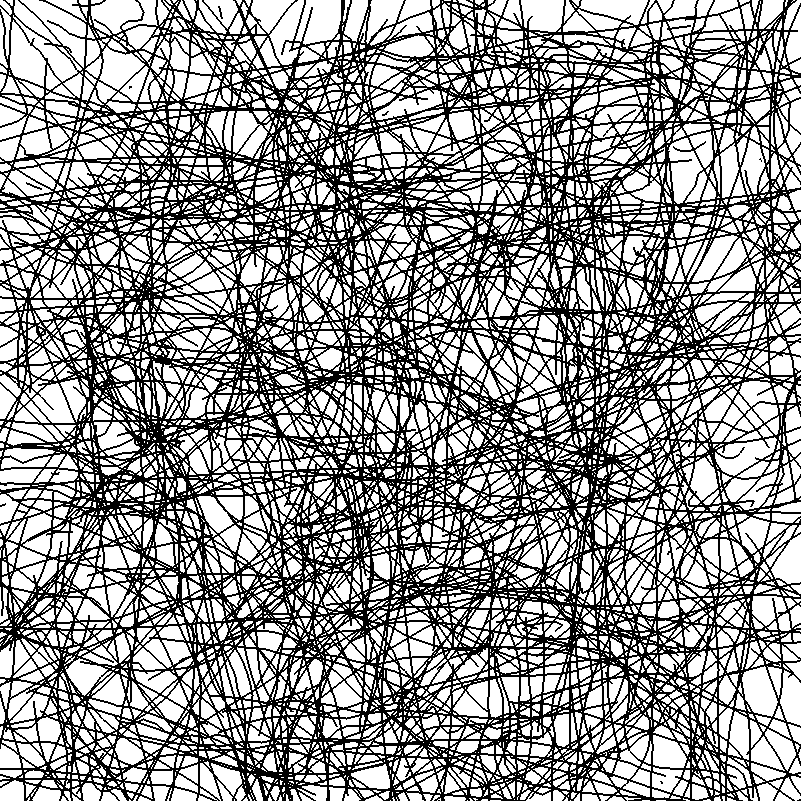}
  \caption{A network model resulting from a paper forming simulation.}
  \label{fig:paper}
\end{figure}

In order to reduce complexity in computer simulation, partial differential equation (PDE) models are sometimes replaced by simpler spatial network models. For instance, instead of modelling the permeability pointwise in a porous media simulation, the effective permeability between subregions can be modeled using a weighted graph, see e.g.~\cite{Chu}. Another example is fiber based materials, like paper and cardboard, where individual fibers can be modeled as one-dimensional objects instead of three-dimensional hollow cylinders, resulting in a spatial network model of edges (fibers) and nodes (connections between fibers), see \cite{BIT}. Still the resulting network models are often large and the data (weights) may vary rapidly in space. In a network model of a porous medium the weights may represent spatially varying permeability and in fiber based materials varying fiber width, see Figure~\ref{fig:paper} for an illustration and \cite{paper} for more details. 

This paper considers spatial network models that arise from applications modeled by elliptic partial differential equations, such as heat conduction and linear elasticity. The resulting spatial networks are assumed to be related to weighted graph Laplacians in the scalar case and weighted graph Laplacians in the coordinate directions in the vector-valued case. For the underlying PDE models, homogenization theory efficiently handles numerical upscaling when the data variation is periodic. For non-periodic data, there are various numerical approaches. Successful numerical algorithms typically use parallelization and discretizations on multiple scales. This is exploited in geometric multigrid \cite{multigrid} and domain decomposition \cite{KY16} but also in numerical homogenization techniques such as the multiscale finite element method \cite{msfem}, gamblets \cite{gamblets}, and localized orthogonal decomposition (LOD) method \cite{aprioriLOD,bookLOD}. It is natural to define coarser scales in PDE problems, at least if the geometry is simple, using nested meshes. In this work, we want to apply numerical homogenization techniques inspired by the PDE community to spatial network models, where it is less obvious how coarse scales can be introduced. 

Numerical homogenization techniques from the PDE community have been applied to spatial network models before. In \cite{Ewing,Iliev} the heat conductivity of a spatial network is studied. In these works, local solutions enable the construction of an effective global thermal conductivity tensor. Traffic flows models in \cite{Rossa} consider a governing PDE for the macroscale by formulating traffic flow equations for single network nodes by interpreting the relations as finite difference approximations. The macroscale parameters are computed using a two-scale averaging technique. In \cite{Chu} spatial network models of flows in a porous medium are studied. The network nodes represent pores and the edges represent throats. The microscale model is based on mass conservation equations for the flow through the network. In \cite{BIT} we consider a specific fiber network model of paper and derive a LOD-based numerical method for efficient numerical simulation. However, the key results needed to prove convergence of the proposed method are left as open problems. In this paper we take advantage of the recent work on  domain decomposition methods for spatial network models, see \cite{GoHeMa22}, in order to prove optimal order convergence of the LOD method applied to spatial network models.
    
We consider a spatial network, defined by a symmetric network matrix $\KK$, for which we want to solve an equation of the form: find $\uu\in\VV$ such that for all $\vv\in\VV$,
$$(\KK \uu,\vv)=(\ff,\vv),$$
given right hand side data $\ff$ and with $(\cdot,\cdot)$ denoting the Euclidean scalar product. We apply the LOD method and introduce an artificial coarse-scale using minimal assumptions on the relation between the coarse-scale representation and the network, see \cite{GoHeMa22}. The fine-scale space is defined as the kernel of a projective quasi-interpolation operator onto the coarse-scale, and the multiscale space is the orthogonal complement to the fine-scale space with respect to the inner product induced by $\KK$. In order to show optimal order convergence, we need to show that the basis spanning the multiscale space decays in space. This decay is possible to establish under mild assumptions on the homogeneity, connectivity, and locality of the network. In order to analyze the error in the proposed method, we prove an interpolation error bound in the spatial network setting. The theoretical findings show how the density variation and connectivity of the network affect the approximation properties of the proposed method. The main result is an optimal order a priori error bound in the norm induced by $\KK$. Finally, the theoretical results are illustrated by three numerical examples.

The paper is organized as follows. Section 2 is devoted to preliminary notation and problem formulation. Section 3 introduces coarse finite element spaces and proves an interpolation error bound. In Section 4, the LOD method is presented and an a priori error bound is derived. Finally, in Section 5 numerical examples are presented.

\section{Problem formulation}
This section presents notation, network operators, function spaces, and finally the model problem considered in this paper, with three examples. 

\subsection{Network and operators}
\begin{figure}
  \centering
  \includegraphics[width = 0.5\textwidth]{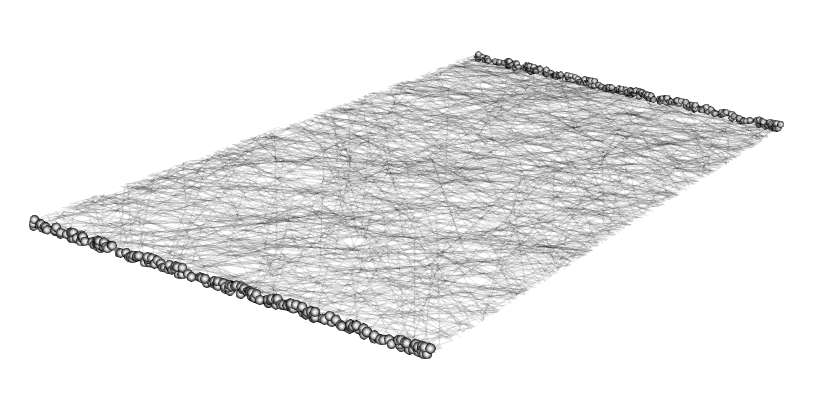}
  \caption{A spatial network with Dirichlet nodes marked on opposite boundary segments.}
  \label{fig:bnd}
\end{figure}

We consider spatial networks represented as connected graphs $\mathcal{G} = (\mathcal{N}, \mathcal{E})$, where the node set
$\mathcal{N}$ is a finite set of points
$x \in \RR^d$ and the edge set
$$\mathcal{E} = \{ \{x, y\} \,:\, \text{an edge connects } x
\text{ and } y \}$$ consists of unordered node pairs. The notation $x \sim y$ means that $\{x, y\}$ is an edge in
$\mathcal{E}$, i.e. $x$ and $y$ are
adjacent. For simplicity we assume that the network reside in the hyper-rectangle $$\Omega = [0,l_1]\times [0,l_2]\times\cdots\times[0,l_d],$$ however, the methodology can be generalized to polygonal and polyhedral domains.  For each pair of adjacent nodes $x\sim y$ we write the Euclidean distance between the nodes as $|x-y|$. Furthermore, we let $\Gamma \subset \partial \Omega$ be the boundary segment where Dirichlet boundary conditions are applied. See Figure~\ref{fig:bnd} for an illustration.

Let the function space $\hat V$ be the space of real valued functions that are defined on the node set $\NN$, and introduce the constrained space 
$$V = \{ v \in \hat V\,:\, v(x) = 0, x \in \Gamma\}.$$
In order to refer to a subset of nodes in the network we define $\mathcal{N}(\omega)=\mathcal{N}\cap\omega$ for any $\omega\subset\Omega$.
Using this notation, let $(u, v)_\omega = \sum_{x \in \mathcal{N}(\omega)} u(x) \cdot v(x)$ and $(u, v) = (u, v)_\Omega$, for all $u,v \in \hat V$.

We further introduce a weighted version of the inner product. This weighted inner product is composed of node-wise diagonal linear operators $M_x: \hat V \to \hat V$:
\begin{equation}
  (M_x v, v) = \frac{1}{2}\sum_{y \sim x} |x - y| v(x)^2.
\end{equation}
These node-wise operators are extended to sets, $\omega \subset \Omega$, by $M_\omega = \sum_{x \in \mathcal{N}(\omega)} M_x$, with $M = M_\Omega$ when the full domain is considered.  The weighted inner product 
$$
(M u,v)=\sum_{x\in\mathcal{N}}(M_x u,v)
$$
defines the norm $|v|_M = (Mv, v)^{1/2}$ and, similarly, $|v|_{M,\omega} = (M_\omega v, v)^{1/2}$ the semi-norms. The squared norm
$|1|_{M,\omega}^2$ of the constant function $1 \in \hat V$ can be interpreted as the mass of the network in subdomain $\omega$.

Next, we define the reciprocal edge-length weighted graph Laplacian. Let
$L_x : \hat V \to \hat V$ be symmetric and defined by
\begin{equation}
  (L_x v, v) = \frac{1}{2}\sum_{y\sim x} \frac{(v(x) - v(y))^2}{|x - y|}
\end{equation}
and introduce the short-hand notation for subdomains
$\omega \subset \Omega$, $L_\omega = \sum_{x \in \mathcal{N}(\omega)} L_x$, and
$L = L_\Omega$. Note that
$(L_\omega u)(x)$ is nonzero for vertices $x$ outside
$\omega$ that are adjacent to nodes in $\omega$.
Further, the weighted Laplacian has the 
constant functions of $\hat V$ in its kernel and defines the semi-norms
$|v|_L = (Lv, v)^{1/2}$ and $|v|_{L,\omega} = (L_\omega
v, v)^{1/2}$.

\subsection{Vector-valued functions}
The models considered in this work are both scalar-valued (e.g., heat conduction) and vector-valued (e.g., structural problems), so we need to introduce vector-valued function spaces and network operators. Let the integer $n$ denote the number of components in the function space
of interest. We introduce 
$$\VV = V^n = V \times \cdots \times V$$ 
as the admissible function space for the unknown and  $\hatVV = {\hat V}^n$ (so that $\VV \subset \hatVV$) as the full space. In
general, the components of $\VV$ can be different by applying 
individual Dirichlet boundary conditions, however,
we assume, for simplicity, that all components are the same. A function
$\vv \in \hatVV$ consists of the components
$v_{1}, v_{2}, \ldots, v_{n}$ and we write $\vv = [v_{1}, v_{2}, \ldots, v_{n}]$. We introduce
$\LL_x : \hat{\VV} \to \hat{\VV}$ as $L_x$ applied component-wise, i.e.
$$\LL_x\vv = [L_x v_{1}, \ldots, L_x v_{n}].$$ 
For $\omega \subset \Omega$, we let
$\LL_\omega = \sum_{x \in \mathcal{N}(\omega)} \LL_x$ and $\LL = \LL_\Omega$. The inner product on this product space is defined in the natural way,
$$(\uu, \vv) = \sum_{i=1}^n(u_{i}, v_{i}),$$
and the semi-norms,
$$|\vv|_{\LL,\omega} = (\LL_\omega \vv, \vv)^{1/2} =
\left(\sum_{i=1}^n|v_i|_{L,\omega}^2\right)^{1/2}$$
where $|\vv|_{\LL} = |\vv|_{\LL,\Omega}$. A similar
notation with subscripts is also used for $\MM$.

\subsection{Model problem} \label{section:problem}
The model problem is expressed with a linear operator
$\KK : \hat{\VV} \to \hat{\VV}$ and a function $\ff \in \hat\VV$: 
\begin{equation}\label{eq:weakform}
  \text{find } \uu \in \VV \,:\, (\KK \uu, \vv) = (\ff, \vv) \text{ for all } \vv \in \VV.
\end{equation}
Since $\uu \in \VV$, $\uu$ is zero on the Dirichlet boundary nodes $\mathcal{N}(\Gamma)$, where
$\Gamma \subset \partial \Omega$. We can easily handle non-zero boundary data 
$\uu(x) = \gg(x)$ for $x \in \mathcal{N}(\Gamma)$ by extending
$\gg$ to all nodes and write $\uu=\uu_0+\gg$, where
$\uu_0 \in \VV$ solves equation (\ref{eq:weakform}) with modified right hand side $(\tilde\ff,\vv):=(\ff,\vv)-(\KK\gg,\vv)$.

Next we make some assumptions on the operator $\KK$.
\begin{assumption}\label{ass:K} The operator $\KK$
  \begin{enumerate}
  \item is bounded and coercive on $\VV$ with respect to $\LL$, i.e.\
    there are constants $0<\alpha \leq \beta < \infty$ such that
    \begin{equation}\label{eq:equivKL}
      \alpha (\LL \vv, \vv) \le (\KK \vv, \vv) \le \beta (\LL \vv, \vv)
    \end{equation}
    for all $\vv \in \VV$, and
  \item can be written as a sum $\KK = \sum_{x \in \mathcal{N} }\KK_x$ of
    operators $\KK_x\,:\, \hat{\VV} \to \hat{\VV}$, where $\KK_x$ are symmetric
    positive semi-definite and only depends and has support on $x$ and nodes adjacent to $x$.
  \item admits a unique solution to equation (\ref{eq:weakform}).
  \end{enumerate}
\end{assumption}
The operator $\KK$ is symmetric as a consequence of $\KK_x$ being
symmetric and the bilinear form $(\KK \cdot,\cdot)$ is an inner product on $\VV$. With
$\KK_\omega = \sum_{x \in \mathcal{N}(\omega)} \KK_x$ we define the following (semi-)norms
$|\vv|_{\KK} = (\KK \vv, \vv)^{1/2}$ and $|\vv|_{\KK, \omega} = (\KK_\omega \vv, \vv)^{1/2}$. The assumptions on $\KK$ are similar to the once made in \cite{GoHeMa22} where an  iterative method for the same problem class is proposed. We now give three examples of system matrices $\KK$ that we consider in this work.

\begin{example}[Heat conductivity] \label{example:heat}
Since this is a scalar example we drop the bold face notation. The same goes for the corresponding (first) numerical example on this model in Section \ref{sec:num}. 
  Let $n=1$ and $u$ be the sought temperature distribution in the nodes of
  the network. We define
  \begin{equation*}
    (K_x v,v) = \frac{1}{2}\sum_{y \sim x} \gamma_{xy} \frac{(v(x)-v(y))^2}{|x - y|},
  \end{equation*}
  where $0 < \gamma_{xy} < \infty$ is heat conductivity on the edges. Assumption~\ref{ass:K} is
  satisfied with $\alpha = \min_{x \sim y} \gamma_{xy}$ and $\beta = \max_{x \sim y} \gamma_{xy}$. By placing at
  least one of the nodes at the Dirichlet boundary $\Gamma$, constants are removed from $V$ and the kernel of $L$ restricted to $V$
  contains only zero and we therefore have a unique solution. The right hand
  side $f$ represent an external heat source.
\end{example}

\begin{example}[Spring model] \label{ex:spring}
  Let $d = n = 3$, and $\partial_{xy} = |x - y|^{-1}(x - y)$ be
  the unit direction vector for edge $\{x, y\}$, then we can define
  \begin{equation} \label{eq:spring}
    (\KK_x\vv,\vv) =\frac{1}{2} \sum_{y \sim x} \gamma_{xy}  \frac{((\vv(x)-\vv(y))^T\partial_{xy})^2}{|x - y|},
  \end{equation}
  where $0 < \gamma_{xy} < \infty$ measures the elasticity of the
  edges. The upper bound of the first assumption in Assumption~\ref{ass:K} is satisfied
  with $\beta = \max_{x \sim y} \gamma_{xy}$ since
  $\partial_{xy}$ has unit length. Whether the lower bound is satisfied
  or not depends on the geometry of the network. At least $d$ nodes
  need to be in $\Gamma$ and they have to span a plane. Additionally,
  the network needs to be a rigid structure. The value of $\alpha$
  depends on $\gamma_{xy}$ and on the structural rigidity of the
  network. We seek the displacement $\uu$ of the nodes under the load $\ff$.
\end{example}

\begin{example}[Fiber network model] \label{ex:fiber}
  Example \ref{ex:spring} can be expanded to represent beams by adding bending stiffness to the edges. A linearized Euler--Bernoulli model can be written on a similar form as \eqref{eq:spring}. For $k=1,2$, then
  \begin{equation} \label{eq:fiber} 
  \begin{split}
    &(\KK^{(k)}_x\vv,\vv) = \\
    &\sum_{\substack{y\sim x \wedge z\sim x \\ y \not=z}} \gamma_{xyz}^{(k)} \frac{|x-y|+|x-z|}{4}  \left(\frac{(\vv(y)-\vv(x))^T\eta_{xyz}^{y,(k)}}{|x-y|}+\frac{(\vv(z)-\vv(x))^T\eta_{xyz}^{z,(k)}}{|x-z|}\right)^2,
  \end{split}
  \end{equation}
   where $\eta_{xyz}^{y,(1)}=\eta_{xyz}^{z,(1)}$ is a unit vector orthogonal to both $\partial_{xy}$ and $\partial_{xz}$, and $\eta_{xyz}^{r,(2)} = \partial_{xr}\times \eta_{xyz}^{r,(1)}$ for $r = y,z$. Adding components $\KK_x^{(1)}$ and $\KK_x^{(2)}$ to \eqref{eq:spring} results in an operator that capture tensile strains in the edges with the spring model and bending resistance with the addition. For more details about this network model, see \cite{BIT}.
\end{example}

\section{Coarse scale representation}
\label{sec:coarse_scale_representation}

The aim of this work is to derive an upscaled representation of the spatial network model problem \eqref{eq:weakform} using the LOD methodology. This representation should have significantly fewer degrees of freedom, but still yield an accurate solution to the original problem.
In this section, using a construction first presented in \cite{GoHeMa22}, we define a coarse scale finite element representation that will be used to form the LOD space. The construction involves three main steps. First, we make assumptions on the spatial network, since not all networks allow for accurate upscaling. In essence, the network should resemble a homogeneous material on the coarse scale. Second, a finite element mesh and coarse function space is introduced on the coarse scale. Third, we introduce a novel idempotent interpolation operator onto the finite element function space and establish the corresponding interpolation error bound.

\subsection{Network assumptions}\label{sec:network}
Four assumptions on the network are made, guaranteeing homogeneity,
connectivity, and locality on a coarse scale. As a technical tool for the assumptions, and later for the definition of the finite element mesh, we define boxes
$B_R(x) \subset \Omega$ with side length $2R$ and midpoint $x = (x_1, \ldots, x_d)$ as follows. Let
$$B_R(x) = [x_1 - R, x_1 + R) \times \cdots \times [x_d - R, x_d
+ R),$$ but if $x_i + R = l_i$, we replace $[x_i - R, x_i + R)$ with
$[x_i - R, x_i + R]$.

From \cite{GoHeMa22} we recall the following network assumptions.
\begin{assumption}[Network assumptions]
  \label{ass:network}
  There is a length-scale $R_0$, a uniformity constant $\sigma$, and a
  density $\rho$, so that
  \begin{enumerate}
  \item (homogeneity) for all $R \ge R_0$ and $x \in \Omega$, it holds
    that
    \begin{equation*}
      \rho \le (2R)^{-d}|1|^2_{M,B_R(x)} \le \sigma \rho,
    \end{equation*}
  \item (locality) the edge length $|x - y| < R_0$ for all edges $\{x, y\} \in \mathcal{E}$,
  \item (boundary density) for any $y \in \Gamma$, there is an
    $x \in \mathcal{N}(\Gamma)$ such that
    $|x - y| < R_0$.
  \item (connectivity) for all $R \ge R_0$ and $x \in \Omega$, there is
    a connected subgraph
    $\mathcal{\bar{G}} = (\mathcal{\bar{N}}, \mathcal{\bar{E}})$ of
    $\mathcal{G}$, that contains
    \begin{enumerate}
    \item all edges with one or both endpoint in $B_R(x)$,
    \item only edges with endpoints contained in $B_{R+R_0}(x)$.
    \end{enumerate}
  \end{enumerate}
\end{assumption}
The four assumptions can be interpreted at scale $R_0$ as follows. The
homogeneity assumption implies that the spatial network has
homogeneous density over the domain in terms of the $M$-norm mass. The
locality assumption says that edges connect only nodes close to each
other, while the boundary density requires that the boundary conditions
are given close enough to nodes. The connectivity assumption
guarantees that nodes close to each other spatially are also close to
each other in the network.

Under the assumptions above, the following Friedrichs and Poincaré inequalities
are proven in \cite{GoHeMa22}.
\begin{lemma}[Friedrichs and Poincaré inequalities]
  \label{lem:poincare}
  If Assumption~\ref{ass:network} holds, then there is a
  $\mu < \infty$ such that for all $R \ge R_0$ and $x \in \Omega$ for
  which
  \begin{itemize}
  \item (Friedrichs) $B_{R}(x)$ contains boundary nodes, it holds that
    \begin{equation*}
      |v|_{M,B_R(x)} \le \mu R |v|_{L,B_{R + R_0}(x)},
    \end{equation*}
    for all $v \in V$,
  \item (Poincaré) $B_{R}(x)$ may or may not contain boundary nodes, it holds that
    \begin{equation*}
      |v - c|_{M,B_R(x)} \le \mu R |v|_{L,B_{R + R_0}(x)},
    \end{equation*}
    for some constant function $c = c(R, x, v)$, for all $v \in \hat V$.
  \end{itemize}
\end{lemma}
The constant $\mu$ enters the interpolation bounds and consequently
affects the accuracy of the homogenization method presented in Section \ref{section:lod}. For simple networks, such
as regular grids, the constant can be shown to be small, while for
most networks, theoretical bounds are generally difficult to
obtain. The constant $\mu$ can, however, be estimated numerically.

As an illustration on how to estimate $\mu$ in the Poincaré case, we
pick an $R$ and an $x$ far from the boundary. Then we study the
eigenvalue problem $\bar L v = \lambda \bar M v$, $v \in \hat V$,
where $\bar M$ and $\bar L$ are the mass operator and reciprocal
edge-length weighted graph Laplacian for the graph $\mathcal{\bar G}$
in Assumption~\ref{ass:network}. We note that
\begin{equation}
  \label{eq:eigenvalue_neumann}
    |v - c|_{M,B_R(x)}^2  \le (\bar M (v - c), v - c)
    \le \lambda_2^{-1} (\bar L v, v) \le \lambda_2^{-1} |v|_{L,B_{R+R_0}(x)}^2
\end{equation}
where $\lambda_2$ is the second smallest eigenvalue, since the
constant $c$ is an eigenvector for the only zero eigenvalue. Thus,
$\mu^2 R^2$ is bounded by the maximum $\lambda_2^{-1}$ attained for
any $R$ and $x$ sufficiently far from the boundary. (A similar
eigenvalue problem can be formulated for the Friedrichs case. See
\cite{GoHeMa22} for details.) Next, we illustrate by a number of
examples of how the connectivity and homogeneity constants can be
estimated numerically based on this eigenvalue problem.

\begin{example}[Numerical estimates of homogeneity and connectivity] \label{ex:eig}
\begin{figure}
  \centering
  \includegraphics[width = 0.3\textwidth]{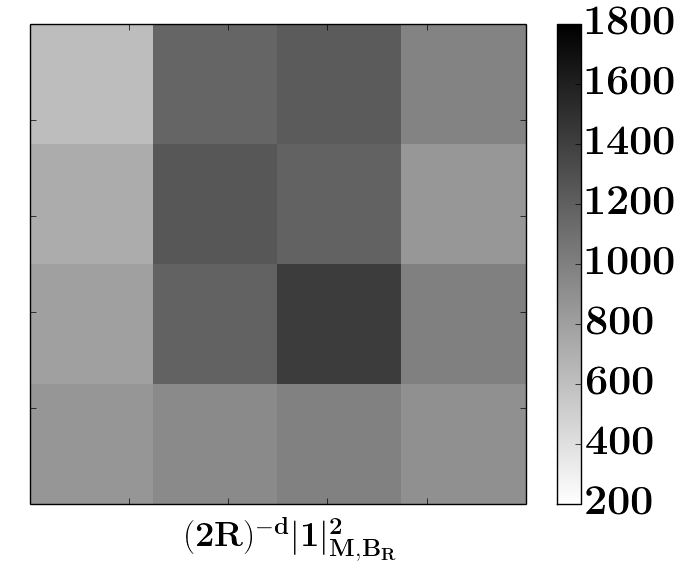}
  \includegraphics[width = 0.3\textwidth]{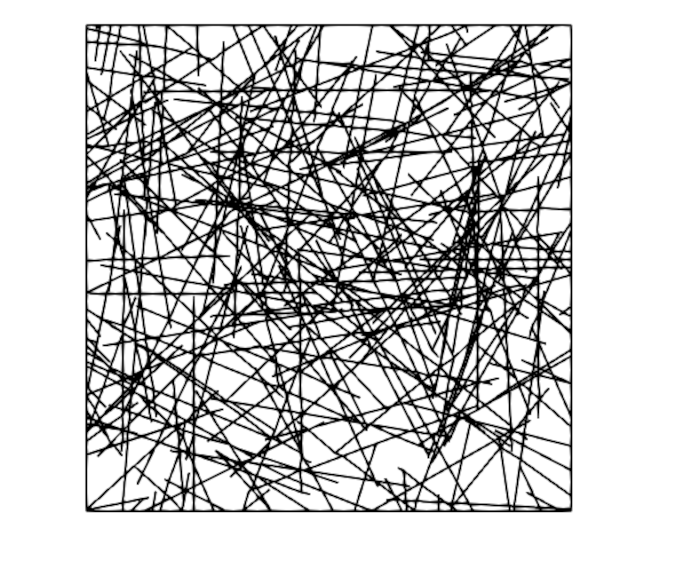}
  \includegraphics[width = 0.3\textwidth]{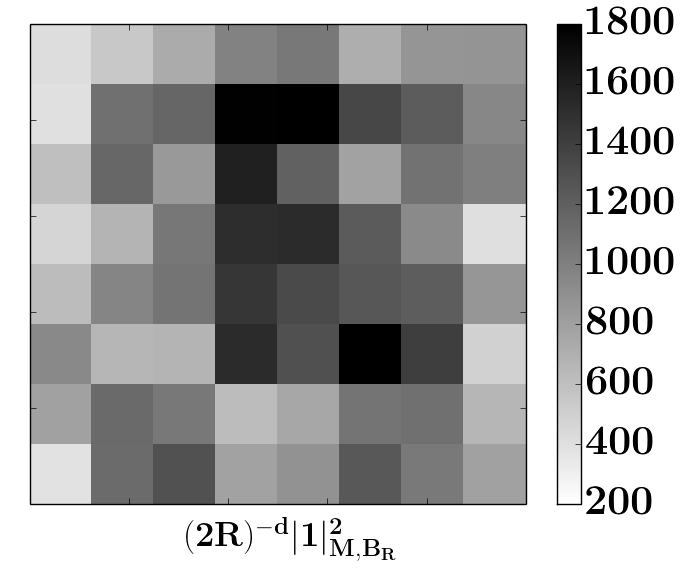}
  \caption{Illustration of the homogeneity equation in Assumption 3.1 for a random network generated on the square $[0,0.1]^2$ with edge length $r=0.05$ and density $\rho_{\text{gen}} = 10^3$ (center) for a grid of $B_R$ with $R=1.25\cdot 10^{-2}$ (left) and $R=6.25\cdot 10^{-3}$ (right).}
  \label{fig:homoex}
\end{figure}

To visualize the homogeneity and connectivity assumptions of Assumption~\ref{ass:network}, several random square networks are generated and evaluated. The analyzed networks are generated by randomly placing edges with a fixed length in a domain. Three attributes categorize each network: the side length $R$ of the domain, the length $r$ of the edges placed, and density $\rho_{\text{gen}}$. The networks are generated in three steps. First, the edges are randomly placed with their midpoints in the extended domain $[-r,R+r]^2$ with a random rotation. This extension guarantees uniform coverage, and any part of an edge placed outside the network domain $[0,R]^2$ is removed. Edges are placed until the total edge length of the network is $\rho_{\text{gen}}R^2$. The second step is adding a node in every point where two edges intersect. The final step removes any loose edges and combines nodes closer than $0.01r$ to guarantee a lower bound on the edge lengths. The largest remaining connected network is kept. An illustration of the homogeneity assumption is shown in Figure~\ref{fig:homoex} for a network with parameter $R=0.1,r=0.05,$ and $\rho_{\text{gen}} = 10^3$. This figure shows how the mean value of $\rho_{\text{gen}}\approx \frac{|1|^2_{M,B_R}}{(2R)^2}$ stays similar when $R$ is halved but varies more throughout the network. The connectivity property is analyzed for multiple networks with multiple parameters, and a composite of the results is presented in Figure~\ref{fig:connex}. We see that $\lambda_2^{-1}$ scales with $R^2$ and thus the connectivity assumption holds for the networks in the interval of $R$ analyzed with $\mu^2 \approx 10$.

\begin{figure}
  \centering
  \includegraphics[width = 0.3\textwidth]{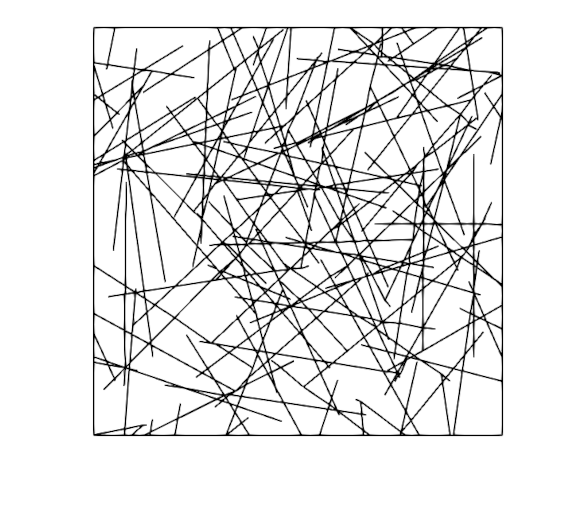}
  \includegraphics[width = 0.3\textwidth, height = 0.28\textwidth]{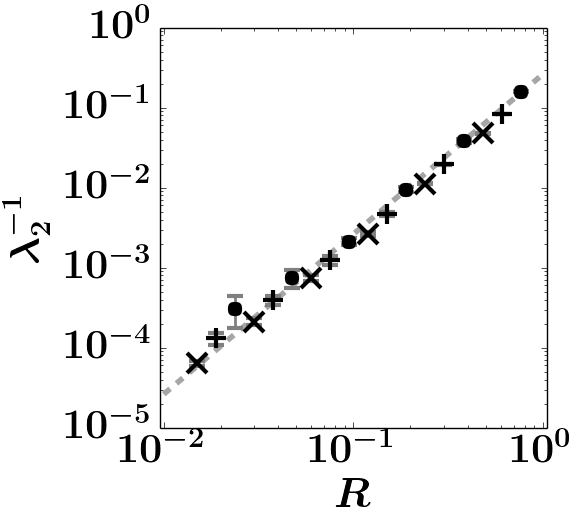}
  \includegraphics[width = 0.3\textwidth]{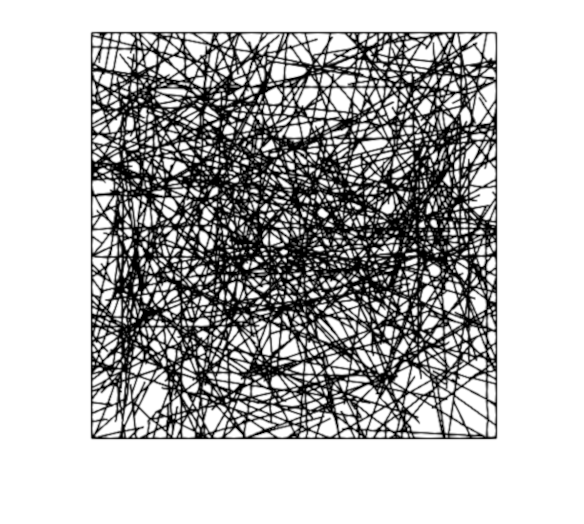}
  \includegraphics[width = 0.5\textwidth]{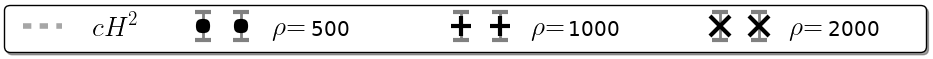}
  
  \caption{The eigenvalue $\lambda_2^{-1}$ for different networks with varying sizes and densities. Each combination is analyzed ten times, with the mean (marker) and standard deviation (feet) results presented. Comparable networks to the one in Figure~\ref{fig:homoex} are shown with $\rho_{\text{gen}} = 5\cdot 10^2$ (left) and $2\cdot 10^3$ (right). }
  \label{fig:connex}
\end{figure}

\end{example}

\subsection{Coarse mesh}\label{sec:mesh}
With the network embedded in a domain $\Omega$ we can introduce a
family of meshes over $\Omega$ for the coarse discretization. The
elements have to be larger than the length-scale $R_0$ of the
network introduced in the previous section. To help
convey the main message of the paper, we choose a simple mesh of
hypercubes (squares for $d = 2$, and cubes for $d = 3$, etc.).  For a general polygonal or polyhedral domain triangles or tetrahedrons can be used. The main difference in the analysis is that constants also will depend on the shape regularity parameter of the mesh.

\begin{figure}
  \centering
  \includegraphics[width = 0.3\textwidth]{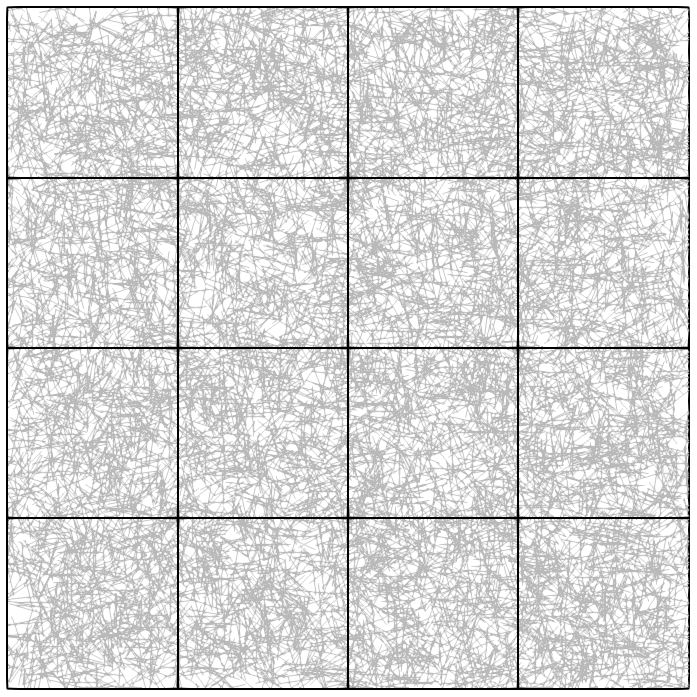} \qquad\qquad \includegraphics[width = 0.3\textwidth]{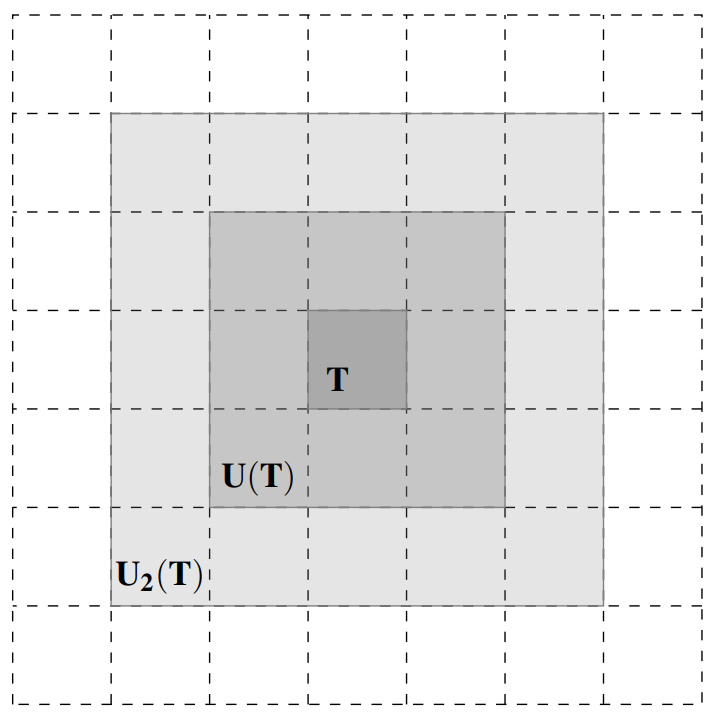}
  \caption{A partition $\mathcal{T}_{H}$, with $H = 1/4$ on the unit square (left) and an illustration showing the recursive operator $U_k$ (right).} 
  \label{fig:partition}
\end{figure}

Let $\mathcal{T}_H$ be subdivisions of
$\Omega$ into elements of side length $H$ as follows,
\begin{equation*}
  \mathcal{T}_{H} = \{ B_{H/2}(x) \,:\, x = (x_1, \ldots, x_d) \in \Omega \text{ and } H^{-1}x_i + 1/2 \in \mathbb{Z} \text{ for } i = 1, \ldots, d\}.
\end{equation*}
We require that $l_1, \ldots, l_d$ are integer multiples of $H$ so that the mesh covers $\Omega$. 
The box definition $B_{R}(x)$ from the previous section is used here. This makes $\mathcal{T}_H$ a true partition so that
each point in $\Omega$ is included in exactly one element. An illustration of such a partition is presented in Figure~\ref{fig:partition} (left). We assume that the boundary segments $\Gamma$ are union of mesh element edges (or faces) so that a conforming finite element function space can be defined.

The elements must be larger than the length-scale $R_0$ of the network. In fact, to define a stable idempotent quasi-interpolation operator, we require that
\begin{equation}
  \label{eq:Hnoll}
  H \ge 4dR_0.
\end{equation}
For a motivation of this lower bound, see the proof of Lemma~\ref{lem:interpolation_stability}. Thus, these meshes are coarse compared to the homogeneity and connectivity length-scale
$R_0$. Note, we do not require that the mesh nodes coincide with the network nodes.

To handle patches of elements in a mesh $\mathcal{T}_H$, we introduce
the notation $U$. We let, for $\omega \subset \Omega$,
$$U(\omega) = \{ x \in \Omega \,:\, \exists T \in \mathcal{T}_H \,:\, x \in T,  T \cap \overline \omega
\ne \emptyset\}.$$ For instance, $U(T)$ contains
the points both in $T$ and in the elements adjacent to $T$. Recursively, we define $U_j(\omega) = U_{j-1}(U(\omega))$ with
$U_1 = U$. An illustration of $T, \ U(T)$, and $U_2(T)$ can be found in Figure~\ref{fig:partition} (right).


\subsection{Interpolation}
In this section, we define the function space to be used for the
coarse representation and an interpolation operator from the functions
on the network to this coarse space. The interpolation operator is of
Scott--Zhang type (see \cite{SZ90}) and is defined by use of an
$M$-dual basis. By showing that the dual basis functions are
appropriately bounded, we can use the results in \cite{GoHeMa22} to
obtain the accuracy and stability result for the interpolation
operator in Lemma~\ref{lem:interpolation}. A Scott--Zhang type
operator is used here instead of Clément operator since it is
important for the analysis of the LOD method that the interpolation
operator is idempotent.

Let $\mathcal{\hat Q}_H$ denote the continuous real functions over
$\Omega$ whose restriction to $T \in \mathcal{T}_H$ can be written as
a linear combination of $z = (z_1, \ldots, z_d) \mapsto z^{\alpha}$
for multi-index $\alpha$ with $\alpha_i \in \{0,1\}$,
$i=1,\dots,d$. For $d = 2$, this is the space of functions that are
bilinear on each element. The functions satisfying the boundary
conditions are
$\mathcal{Q}_{H} = \{ p \in \mathcal{\hat Q}_H \,:\, p|_\Gamma = 0
\}$. We let $\hat V_H$ and $V_{H}$ be the restriction of
$\mathcal{\hat Q}_H$ and $\mathcal{Q}_H$ to the nodes in the network.

\begin{figure}
\centering
\includegraphics[width = 0.4\textwidth]{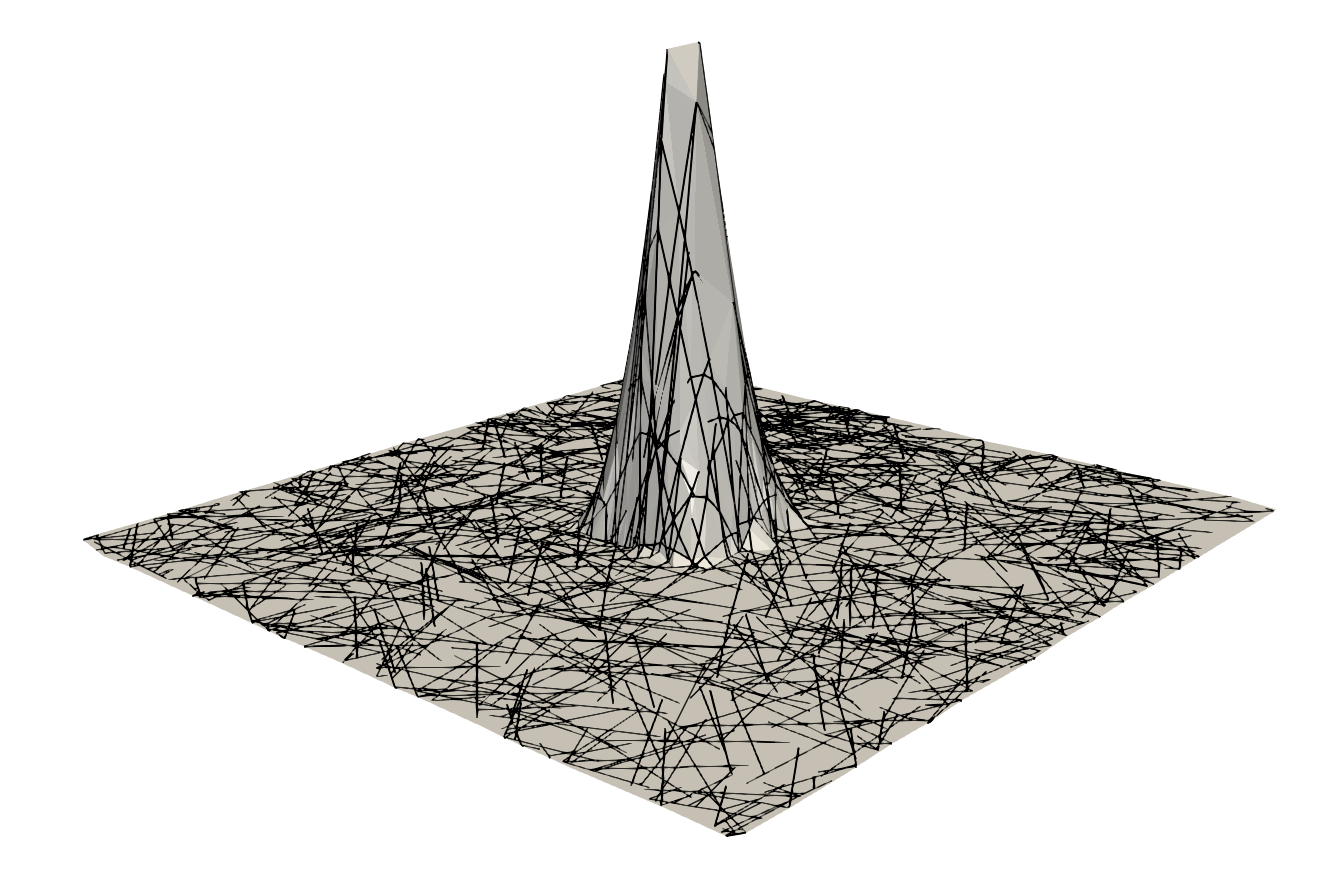}
\caption{A two-dimensional network (black) displaced to the discrete function values of $\phi$ representing the $\varphi$ basis function. The gray shading is a triangulation of the displaced network for illustrative purposes.}
\label{fig:bilinear}
\end{figure}

From this point on, we study a fixed $H$. Denote by
$\phi_1, \ldots, \phi_m \in \mathcal{\hat Q}_{H}$ the Lagrange finite
element basis functions, $\varphi_1, \ldots, \varphi_m \in \hat V_{H}$
their restrictions to the network nodes and $y_1, \ldots, y_m$ the
corresponding mesh nodes. An illustration of a discrete $\varphi_i$
function can be found in Figure~\ref{fig:bilinear}.  We assume that
the basis functions are sorted so that the basis functions
$\phi_1, \ldots, \phi_{m_0}$, $m_0 < m$ span $\mathcal{Q}_{H}$ that
vanish on $\Gamma$.
For each mesh node $y_k$, we denote by $T_k$ the
unique element that contains it and define $\psi_k \in \hat V_H(T_k)$
that satisfies
\begin{equation}
  (M_{T_k} \psi_k, \varphi_\ell) = \delta_{k\ell}
\end{equation}
for all $\ell = 1, \ldots, m$. The interpolation operator is
then defined by
$$\II v = \sum_{k=1}^{m_0} (M_{T_k} \psi_k, v) \varphi_k.$$

We let the subscript $z$ of a constant $C_z$ indicate a dependency on a constant $z$, where the exact value of $C_z$ may differ (by some generic constant). 
\begin{lemma}[Dual basis norm bound]
  \label{lem:interpolation_stability}
  If Assumption~\ref{ass:network} hold and $H\geq4dR_0$, then 
  \begin{equation}
    |\psi_k|_{M,T_k} \le C_d \rho^{-1/2} H^{-d/2}
  \end{equation}
  for mesh nodes $k = 1,\ldots,m$.  
\end{lemma}
\begin{proof}
  To relieve the notation, we omit subscript $k$ and set
  $\psi = \psi_k$ and $T = T_k$ in this proof. Denote the mesh
  vertices in $T$ by $y_1,\ldots,y_{2^d}$ and let $y_1$ be the vertex
  for which $\psi$ is the dual -basis.  We define the positive
  semi-definite Gram matrix $\Lambda$ with entries
  $\Lambda_{ij} = (M_T \varphi_j, \varphi_i)$ for
  $i,j=1,\ldots,2^d$. Let
  $\alpha = (\alpha_1, \ldots, \alpha_{2^d})^T$ and express the
  dual basis as $\psi = \sum_{\ell=1}^{2^d} \alpha_\ell
  \varphi_\ell$. Then by the definition of $\psi$, we have
  $\alpha = \Lambda^{-1} (1, 0, \ldots, 0)^T$ and that the sought
  squared norm
  $(M_T \psi, \psi) = \alpha_1 \le \lambda_1(\Lambda)^{-1}$, where
  $\lambda_1(\Lambda)$ is the smallest eigenvalue of $\Lambda$. To be
  able to bound this eigenvalue from below, we split $\Lambda$ into
  the significant contributions from network nodes close to the
  corners of the element.

  Let
  $\hat T_\ell = \{ x \in T \,:\, y_\ell + (x - y_\ell)/r \in
  \overline{T}\}$ be the points in $T$ that is in an $r$ scaling of
  $\overline T$ in the corner of node $y_\ell$. We set $r = 1/(4d)$, but
  keep writing $r$ for brevity. Since the closure of $\hat T_\ell$
  is a scaling of $\overline T$, it is a hypercube of side length
  $r H \ge R_0$. We set
  $\hat T_0 = T \setminus \hat T_1 \setminus \cdots \setminus \hat
  T_\ell$ and define the symmetric positive semi-definite matrices
  $\Lambda^\ell$ for $\ell = 0, 1, \ldots, 2^d$ with entries
  $\Lambda_{ij}^\ell = (M_{\hat T_\ell} \varphi_j, \varphi_i)$. We can
  then write
  $\Lambda = \Lambda^0 + \Lambda^1 + \ldots + \Lambda^{2^d}$.  For
  brevity, let $w_x = |1|_{M, \{x\}}^2$ for $x \in \NN$. Since
  $x \mapsto \phi_j(x)\phi_i(x)$ is continuous and $\hat T_\ell$ is
  path-connected, by the intermediate value theorem there is an
  $x_\ell \in \hat T_\ell$ such that
  \begin{equation*}
    \Lambda_{ij}^\ell = (M_{\hat T_\ell} \varphi_j, \varphi_i) = \sum_{x \in \NN(\hat T_\ell)} w_x \phi_i(x) \phi_j(x) = \sum_{x \in \NN(\hat T_\ell)} w_x \phi_i(x_\ell) \phi_j(x_\ell).
  \end{equation*}
  Using the properties of the smallest eigenvalues of symmetric real
  operators $A$ and $B$:
  (i)~$\lambda_1(A + B) \ge \lambda_1(A) + \lambda_1(B)$ and
  (ii)~$\lambda_1(\alpha A + \beta B) \ge \min(\alpha,
  \beta)\lambda_1(A + B)$, and that we from
  Assumption~\ref{ass:network}.1 have
  $\sum_{x \in \NN(\hat T_\ell)} w_x = |1|_{M, \hat T_\ell}^2 \ge
  r^{d}\rho H^{d}$, we get
  \begin{equation*}
    \lambda_1(\Lambda) \ge \lambda_1(\Lambda^0) + \lambda_1(\Lambda^1 + \ldots + \Lambda^{2^d}) \ge \lambda_1(\Lambda^1 + \ldots + \Lambda^{2^d}) \ge r^{d}\rho H^{d} \lambda_1(G)
  \end{equation*}
  where $G$ is a matrix with entries
  $G_{ij} = \sum_{\ell = 1}^{2^d} \phi_i(x_\ell) \phi_j(x_\ell)$. The
  next step is to bound $\lambda_1(G)$ from below by means of the
  Gershgorin circle theorem.
  
  We study the first row of $G$ and note that all entries in the row
  are positive. The distance $D_1$ between zero and the Gershgorin disc for
  the first row can be expressed as the difference between
  the diagonal entry and the sum of the (all positive) non-diagonal
  entries on the row, i.e.\
  \begin{equation*}
    \begin{aligned}
      D_1 & = G_{11} - \sum_{\ell = 2}^{2^d} G_{1\ell} = \sum_{\ell=1}^{2^d} \phi_1(x_\ell)^2 - \phi_1(x_\ell) \sum_{j=2}^{2^d} \phi_j(x_\ell) \\
      &= \sum_{\ell=1}^{2^d} \phi_1(x_\ell) (2\phi_1(x_\ell)-1) = \phi_1(x_1) (2\phi_1(x_1)-1) + \sum_{\ell=2}^{2^d} \phi_1(x_\ell) (2\phi_1(x_\ell)-1), \\
    \end{aligned}
  \end{equation*}
  where the partition of unit of the basis functions was used.  Since
  $x_\ell \in \hat T_\ell$, the values of the basis function $\phi_1$
  in these points are bounded as follows
  \begin{equation*}
    \begin{aligned}
      (1 - r)^d \le {}& \phi_1(x_1) \le 1 && \text{and}\\
      0 \le {}& \phi_1(x_\ell) \le r^k && \text{if } y_1 \text{ and } y_\ell \text{ differ in } 1 \le k \le d \text{ components.}
    \end{aligned} 
  \end{equation*}
  The condition for the second bound can also be phrased as $k$ being
  the minimum number of edges of the hypercube $T_\ell$ to traverse to
  reach $y_\ell$ from $y_1$. We note that, for each $k$, there are
  $d \choose k$ element corners $y_\ell$ for which this condition
  hold. This allows us to write
  \begin{equation*}
    \sum_{\ell=2}^{2^d} \phi_1(x_\ell) \le \sum_{k=1}^{d} {d \choose k} r^k = (1 + r)^d - 1,
  \end{equation*}
  which will be useful next.
  
  With the particular choice $r = 1/(4d)$, we use Bernoulli's
  inequality $(1-r)^d \ge 1 - rd = 3/4$ to bound
  $2\phi_1(x_1) - 1 \ge 2(1-r)^d - 1 \ge 1/2$. Using this inequality
  again, together with $2\phi_1(x_\ell)-1 \ge -1$ and
  $(1+r)^d \le e^{rd} = e^{1/4}$, we bound $D_1$ as follows,
  \begin{equation*}
      D_1  \ge \frac{1}{2} \phi_1(x_1) - \sum_{\ell=2}^{2^d} \phi_1(x_\ell) 
       \ge \frac{1}{2} (1-r)^d  - (1 + r)^d + 1 
       \ge \frac{3}{8} - e^{1/4} + 1 > 0.
  \end{equation*}
  Thus, the distance between zero and the Gershgorin disc for the
  first row is bounded below by a positive constant. The argument can
  be repeated for the $2^d$ rows and we get that all eigenvalues of
  $G$ are bounded below, and in particular that $\lambda_1(G) \ge
  C$. We obtain the asserted inequality
  \begin{equation*}
    |\psi|_{M, T}^2 = (M_T \psi, \psi)\le \lambda_1(\Lambda)^{-1} \le (4d)^{d}\rho^{-1} H^{-d} \lambda_1(G)^{-1} \le C_d \rho^{-1} H^{-d}.
  \end{equation*}
\end{proof}

With the above result on the dual basis norm bound, we can use the
interpolation bound for spatial network models established in
\cite{GoHeMa22}.

\begin{lemma} 
  \label{lem:interpolation}
  If Assumption~\ref{ass:network} holds and $H\geq 4dR_0$, then for $v \in V$ and all
  $T \in \mathcal{T}_H$,
  \begin{equation}
    H^{-1}| v - \II v |_{M,T} + |\II v|_{L,T} \le C_{d, \mu, \sigma} |v|_{L,U_3(T)}.
  \end{equation}  
\end{lemma}
This is an element local version of \cite[Lemma~5.4]{GoHeMa22} with a
different choice of nodal variable $v \mapsto (M_{T_k} \psi_k, v)$ for
the interpolation operator. In \cite{GoHeMa22}, a Clément
interpolation operator is used, while a Scott--Zhang interpolation
operator is used here. The proof from \cite{GoHeMa22} can be used
almost verbatim to prove this element local version. We just leave out
the summation over all elements in the end and use bound on the norm
of $\psi_k$ shown in Lemma~\ref{lem:interpolation_stability}.

\section{Numerical homogenization} \label{section:lod}

Given the spatial network model and a coarse scale finite element space, the aim is to derive an accurate upscaled representation of the model problem. This is accomplished by using the localized orthogonal decomposition (LOD) technique originally developed for numerical homogenization of elliptic partial differential equations with heterogeneous data, see \cite{aprioriLOD,bookLOD}. An accurate representation is achieved by decoupling the fine scale computations into local subproblems and thereby constructing a multiscale basis that captures the data variation. The heterogeneities present in the spatial network setting comes from the geometry of the graph and the spatially varying weights. With the results from Section~\ref{sec:coarse_scale_representation}, we can derive the LOD method for the model problem defined in equation (\ref{eq:weakform}). We let $\III \,:\, \VV \to \VV$ be defined as $\III \vv = [\II v_1, \ldots, \II v_n]$ and introduce a fine scale space
\begin{equation*}
  \WW = \ker(\III)=\{\vv\in \VV \,:\, \III \vv = 0\}.
\end{equation*}
\subsection{Ideal multiscale method}
The multiscale space $\VV_{H}^\text{ms}$ is defined as the
orthogonal complement of $\WW$ with respect to the inner product
induced by $\KK$.
For every $\vv \in \VV$ we define a fine scale projection operator $\QQ \,:\, \VV \to \WW$ such that
\begin{equation} \label{eq:QQ}
  (\KK \QQ \vv, \ww) = (\KK \vv, \ww)
\end{equation}
for all $\ww \in \WW$.

\begin{definition}
The ideal multiscale space is defined as 
\begin{equation*}
  \VV_{H}^\text{ms}=\{(1-\QQ)\vv\,:\,\vv\in \VV\}.
\end{equation*}
Any vector $\vv\in \VV$ can be decomposed into
\begin{equation*}
  \vv=(\vv-\QQ \vv)+\QQ \vv\in \VV_{H}^\text{ms}\oplus \WW
\end{equation*}
with the two terms being $\KK$-orthogonal. The ideal multiscale approximation $\uu_H$ of $\uu$ fulfills: 
find $\uu_H\in \VV_{H}^\text{ms}$ such that for all $\vv \in \VV_{H}^\text{ms}$
\begin{equation}\label{eq:u_H}
(\KK \uu_H, \vv) = ( \ff, \vv).
\end{equation}
\end{definition}

\begin{lemma}\label{lem:ideal}
  The error in the approximate solution $\uu_H$, defined in equation (\ref{eq:u_H}), fulfills
  $$
  |\uu-\uu_H|_\KK\leq  C_{\alpha,d,\mu,\sigma} H |\ff|_{\MM^{-1}},
  $$
  where $|\ff|^2_{\MM^{-1}}=(\MM^{-1}\ff,\ff)$.
  \end{lemma}
  \begin{proof}
  The error $\uu-\uu_H\in \WW $ is bounded by
  \begin{align*}
  |\uu-\uu_H|^2_\KK&=(\KK \uu,\uu-\uu_H)\\
  &=(\ff,\uu-\uu_H)\\
  &\leq |\ff|_{\MM^{-1}} |\uu-\uu_H-\III(\uu-\uu_H)|_\MM\\
  &\leq C_{\alpha,d,\mu,\sigma} H |\ff|_{\MM^{-1}} |\uu-\uu_H|_\KK,
  \end{align*}
  where Lemma~\ref{lem:interpolation} is used in all coordinate directions and the overlap of subregions $U_3(T)$ are hidden in $C_{\alpha,d,\mu,\sigma}$. The lemma follows after division by $|\uu-\uu_H|_\KK$.
\end{proof}

For the method to be computationally feasible we need to localize the fine scale correction operators. To do this, we first decompose the computation of $\QQ$ to the elements
$T \in \mathcal{T}_H$. We define $\QQ_T : \VV \rightarrow \WW$, such that for any $\vv\in\VV$ 
\begin{equation}
  (\KK \QQ_T \vv, \ww) = (\KK_T \vv, \ww)
\end{equation}
for all $\ww \in \WW$. By this choice, $\QQ_T \vv$ is independent of the values of $\vv$ in points not
adjacent to nodes in $T$. Note that $\QQ = \sum_{T\in \mathcal{T}_H} \QQ_T$ since $\KK_T=\sum_{x\in T}\KK_x$ sums up to
$\KK$.  

\subsection{The LOD method}
With the fine scale correction decomposed into element components, we want to localize the computation of those components to element patches. Let
\begin{equation*}
  \WW(\omega)=\{\ww\in \WW\,:\,\ww(x)=0\text{ for }x \in \mathcal{N}(\Omega \setminus \omega)\}
\end{equation*}
and let $\QQ^k_T \vv \in \WW(U_k(T))$ be the solution to
\begin{equation}
  (\KK \QQ^k_T \vv, \ww) = (\KK_T \vv, \ww)
\end{equation}
for all $\ww \in \WW(U_k(T))$.  We sum the contributions over the elements
to get the full truncated fine scale projection operator
$\QQ^k \vv =\sum_{T \in \mathcal{T}_H} \QQ^k_T \vv$.
By this construction, $\QQ^k$ is an approximation of
$\QQ$ computed on element patches.

To define the truncated LOD space, we define coarse basis functions
for the free nodes. Let $\vvarphi_j$ for $j = 1, 2, \ldots, nm_0$
enumerate the basis functions span
the range of $\III$. We define them as
$\vvarphi_1 = [\varphi_1, 0, \ldots, 0]$,
$\vvarphi_2 = [\varphi_2, 0, \ldots, 0], \ldots, \vvarphi_{m_0} =
[\varphi_{m_0}, 0, \ldots, 0]$,
$\vvarphi_{m_0+1} = [0, \varphi_1, 0, \ldots, 0]$, etc.\ to
$\vvarphi_{nm_0} = [0, \ldots, 0, \varphi_{m_0}]$.

\begin{definition}
The truncated LOD space is given by
\begin{equation*}
  \VV_{H}^{\text{ms},k}=\text{span}\left(\left\{\vvarphi_j - \sum_{T\in\mathcal{T}_H} \QQ^k_T \vvarphi_j\,:\,j = 1, \ldots, nm_0\right\}\right)
\end{equation*}
and the LOD approximation by: find  $\uu^k_H\in \VV_{H}^{\text{ms},k}$ such that for all $\vv \in \VV_{H}^{\text{ms},k}$,
\begin{equation}\label{eq:u_Hk}
  (\KK \uu^k_H, \vv) = (\ff, \vv) .
\end{equation} \label{eq:loclod}
\end{definition}

\begin{figure}
\centering
\includegraphics[width = 0.4\textwidth]{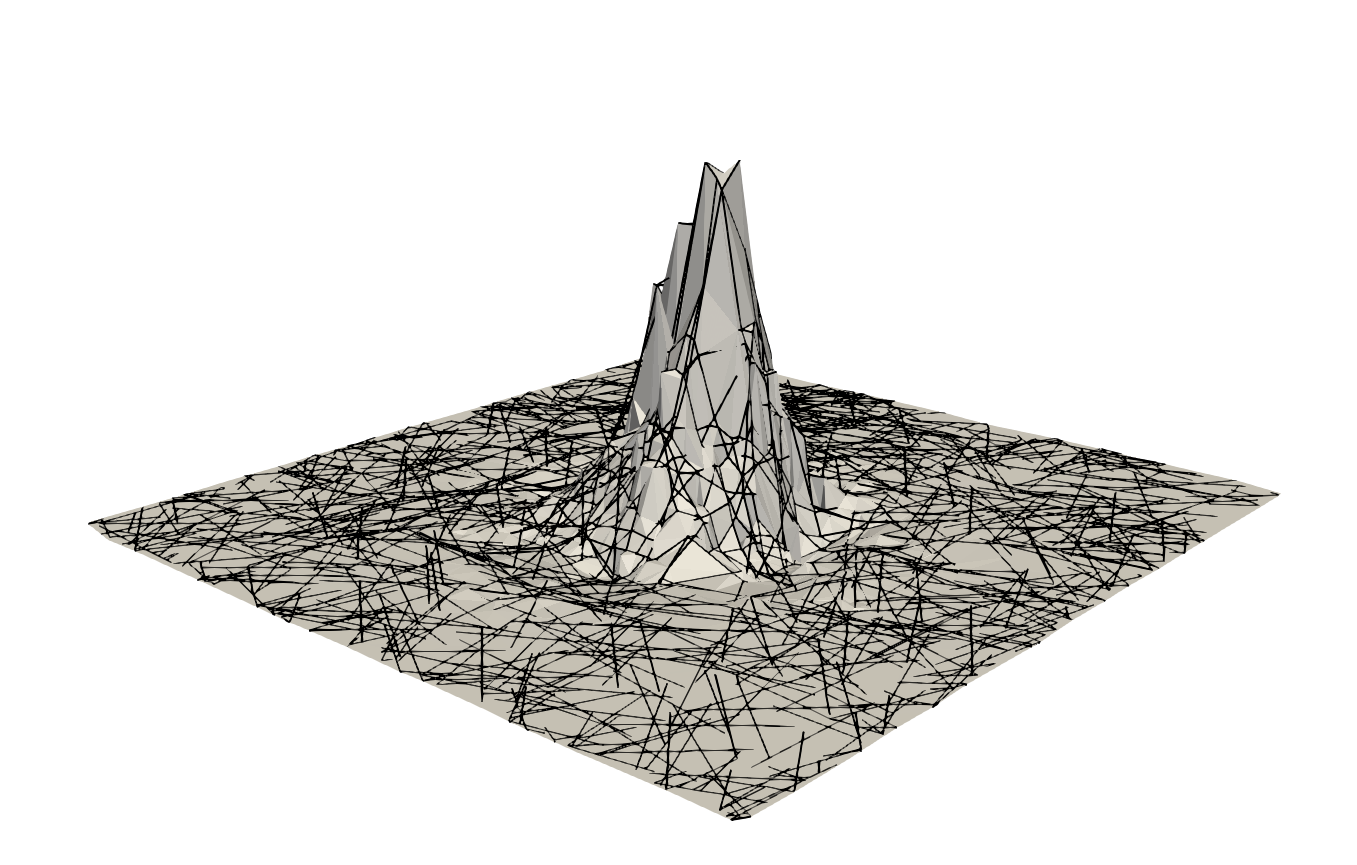}
\hspace{0.05\textwidth}
\includegraphics[width = 0.4\textwidth]{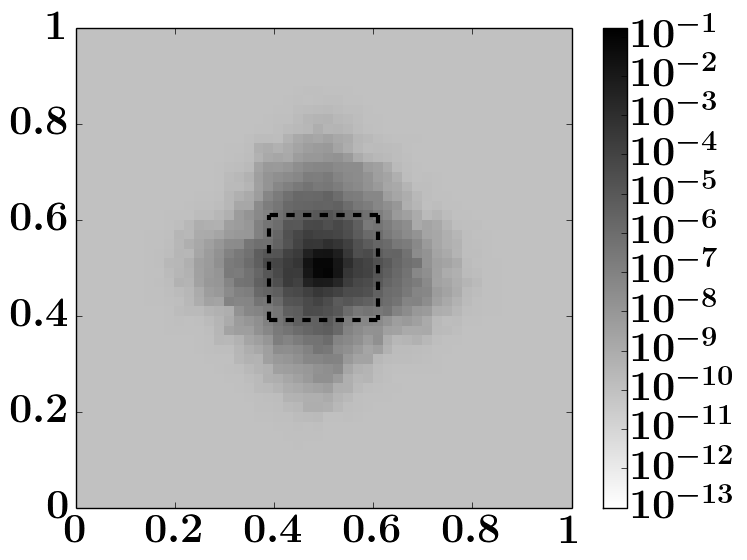}
\caption{The bilinear function, $\varphi$, in Figure~\ref{fig:bilinear} with an ideal fine scale correction, i.e. $(1-\QQ)\vvarphi$ for a heat conductivity problem. The left figure is comparable to Figure~\ref{fig:bilinear}, and the right shows  $(1-\QQ)\vvarphi$ in the entire domain with the area of the left marked with the dashed lines.}
\label{fig:lodbasis}
\end{figure}

\subsection{Decay of fine scale correctors}\label{sec:KY}

In order for $\uu^k_H$ to be a good approximation of $\uu$ for small
values of $k$ we need to show that $\QQ_T \vvarphi_j$ decays quickly away from $T$. This is done by following ideas presented in  \cite{KPY18} and \cite{GoHeMa22}, using Lemma~\ref{lem:interpolation} and a discrete analogue to the product rule.

The idea is to approximate the fine scale projection $\QQ$ using an
iterative domain decomposition technique that spread information
locally in each iteration. By proving that the method converges
quickly we can also draw conclusions about the decay of the $\QQ_T \vvarphi_j$. 


For points $x \in \Omega$, let $U(x)$ be used as short-hand notation
for $U(\{x\})$. In particular, if $y_k$ is a mesh node then $U(y_k)$
is the node patch.
We let $\VV_j = \VV(U(y_j))$ be the space of functions that vanish outside the
node patch for mesh nodes $j = 1, \ldots, m$ (including the fixed mesh nodes $m_0 + 1, m_0 + 2\ldots, m$). The fine scale space $\WW$ is decomposed into
overlapping subspaces
\begin{equation*}
  \WW_j=(1-\III)\VV_j = \{ \vv - \III \vv \,:\,\vv \in \VV_j\}.
\end{equation*}
The relation $\WW_j\subset \WW$ holds since $\III$ is idempotent.
Since the scalar basis $\{\varphi_j\}_{j = 1, \ldots, m}$ is a
partition of unity on $\mathcal{N}$, any $\ww\in \WW$ can be written
as
\begin{equation*}
  \ww=\sum_{j=1}^m(1-\III)(\varphi_j \ww)\quad\text{with}
  \quad(1-\III)(\varphi_j \ww)\in \WW_j,
\end{equation*}
where $\varphi_j \ww = [\varphi_j w_{1}, \ldots, \varphi_j w_{n}]$. 
Remember that $(\III\vv)(y_j)$ is computed by taking a weighted average of $\vv$ in element $T_j$ with $y_j\in T_j$. Therefore, $\III\vv$ has a slightly larger support than $\vv$.  More precisely, any $\ww_j\in\WW_j$ fullfills 
\begin{equation}\label{eq:spread}
\ww_j\in \VV(U_1(\hat T_j))\cap \WW, 
\end{equation}
where $\hat T_j$ is an element adjacent to the node $y_j$.

Now let $\PP_j\,:\,\VV\rightarrow \WW_j$ define the projection such that for any $\vv \in \VV$ and all  $\ww \in \WW_j$
\begin{equation*}
  (\KK \PP_j \vv, \ww) = (\KK \vv, \ww).
\end{equation*}
The operator $\PP=\sum_{j=1}^m \PP_j$ is a preconditioner for $\QQ$, and it is important that it only spreads information a few layers of coarse elements in each application by equation (\ref{eq:spread}). 

Next, we investigate some further properties of $\PP$. The following discrete analogue of a product rule (see \cite{GoHeMa22} for a proof),
\begin{equation}
  \label{eq:product_rule}
  | v \varphi_k |_{L,T}^2 \le 2 \left(H^{-2} | v |_{M,T}^2 + | v |_{L,T}^2\right), \quad k = 1,\ldots,m,\quad T\in\mathcal{T}_H,
\end{equation}
is used in the proofs below.

\begin{lemma}\label{lem:K1K2}
Every decomposition $\ww=\sum_{j=1}^m \ww_j$ with $\ww_j\in \WW_j$ satisfies
\begin{equation*}
  |\ww|_{\KK}^2 \leq C_2 \sum_{j=1}^m |\ww_j|_{\KK}^2
\end{equation*}
and the particular decomposition $\ww_j=(1-\III)(\varphi_j \ww)$ satisfies
\begin{equation*}
  \sum_{j=1}^m |\ww_j|_{\KK}^2 \leq C_1 |\ww|_{\KK}^2.
\end{equation*}
\end{lemma}
\begin{proof}
  For both inequalities, we show them locally first in the $\LL$-norm.
  We start with the first inequality. Pick a $T \in \mathcal{T}_H$. Since
  $\ww_j \in \VV(U_2(y_j))$ and $\LL_T \vv = 0$ for
  $\vv \in \VV(\Omega \setminus U(T))$, we have that $\LL_T \ww_j$ can
  be non-zero for at most $C_d$ mesh nodes $j$, where $C_d$ depends
  only on $d$. Since $\LL_T$ is locally defined in this sense, we get
  \begin{equation*}
    |\ww|_{\LL,T}^2 \le C_d \sum_{j=1}^m |\ww_j|_{\LL,T}^2.
  \end{equation*}
  Summation over $T \in \mathcal{T}_H$ proves the inequality in
  $\LL$-norm. From Assumption~\ref{ass:K} the global $\LL$- and
  $\KK$-norms are equivalent with constant
  $\beta\alpha^{-1}$, thus we get the asserted inequality with
  $C_2 = C_{\alpha, \beta, d}$.

  For the second inequality, we use Lemma~\ref{lem:interpolation}
  componentwise and globally,
  inequality \eqref{eq:product_rule}, the fact that $\III \ww = 0$,
  Lemma~\ref{lem:interpolation} again, and finally a similar locality
  argument of $\LL_T$ as in the previous paragraph and get
  \begin{equation*}
    \begin{aligned}
      \sum_{j=1}^m |\ww_j|_{\LL}^2 & \le C_{d,\mu,\sigma} \sum_{j=1}^m |\varphi_j \ww|_{\LL}^2 \\
      & \le C_{d,\mu,\sigma} \sum_{j=1}^m \sum_{T \in \mathcal{T}_H} \left(H^{-2}|\ww|_{\MM,T}^2 + |\ww|_{\LL,T}^2\right) \\
      & \le C_{d,\mu,\sigma} \sum_{j=1}^m \sum_{T \in \mathcal{T}_H} |\ww|_{\LL,U_3(T)}^2 \\
      & \le C_{d, \mu,\sigma} |\ww|_{\LL}^2.
    \end{aligned}
  \end{equation*}
  We use equivalence of $\LL$- and $\KK$-norms to get the
  second inequality with $C_1 = C_{\alpha,\beta,d,\mu,\sigma}$.
\end{proof}

Using Lemma~\ref{lem:K1K2} one can show the following norm equivalence, where we refer to Lemma~3.1 in \cite{KY16} for a proof of the first statement and the appendix of \cite{GoHeMa22} for a proof of the second statement.
\begin{lemma}\label{lem:normeq}
The following norm equivalence holds
\begin{equation*}
  C_1^{-1}|\ww|_{\KK}^2 \leq (\KK \PP \ww, \ww) \leq C_2 |\ww|_{\KK}^2
\end{equation*}
for all $\ww \in \WW$. Furthermore,
with $\nu=(C_2+C_1^{-1})^{-1}$ and $\ww \in \WW$, it holds
\begin{equation*}
  \sup_{\ww \in \WW} \frac{|(1 - \nu \PP) \ww|_{\KK}}{|\ww|_{\KK}}\leq \gamma<1,
\end{equation*}
where $\gamma\leq \frac{C_2}{C_2+C_1^{-1}}$.
 \end{lemma}

We now define an approximation $\Rdd_T^k\,:\, \VV \rightarrow \WW$ to $\QQ_T$ by the iteration,
\begin{equation}
\Rdd_T^k \vv = \Rdd_T^{k-1} \vv + \nu \PP (\QQ_T - \Rdd_T^{k-1}) \vv,\quad k\geq 1,
\end{equation}
with $\Rdd_T^0 = 0$ and a relaxation parameter $\nu>0$. First we note that $\Rdd_T^k\vv$ is computable without explicitly forming $\QQ_T\vv$ since $\PP_j\QQ_T\vv\in \WW_j$ solves
$$
(\KK\PP_j\QQ_T\vv,\ww_j)=(\KK\QQ_T\vv,\ww_j)=(\KK_T\vv,\ww_j)
$$
for all $\ww_j\in \WW_j$. We further conclude that $\Rdd_T^k$ is local. The right hand side 
$\KK_T\vv$ has support on $U_1(T)$. Since functions in $\WW_j$ have support on $U_1(\hat T_j)$ according to equation (\ref{eq:spread}) only a few of the corresponding projections $\PP_j$ will be non-zero. More precisely $\Rdd_T^1$ will have support on $U_3(T)$ and in general 
$$
\text{supp}(\Rdd_T^k)\subset U_{3k}(T).
$$
We will use this property when we show that $\QQ_T$ decays exponentially. The approximation $\Rdd_T^k$ of $\QQ_T$ fullfills the error bound
\begin{equation}\label{eq:QQdd}
(\QQ_T-\Rdd_T^k)\vv=(1-\nu \PP)(\QQ_T - \Rdd_T^{k-1})\vv = (1-\nu \PP)^k \QQ_T \vv.
\end{equation}

Altogether we get the following approximation result.
\begin{lemma}\label{lem:decay}
It holds
\begin{equation*}
  |(\QQ_T - \Rdd_T^k)\ww|_{\KK} \leq \exp(-k (2C_1C_2)^{-1}) |\ww|_{\KK,T}.
\end{equation*}
\end{lemma}
\begin{proof}
Using equation (\ref{eq:QQdd}) and Lemma~\ref{lem:normeq} we conclude
\begin{equation*}
|(\QQ_T - \Rdd_T^k)\ww|_{\KK} \le \gamma^k |\QQ_T \ww|_{\KK} \le \gamma^k |\ww|_{\KK,T}.
\end{equation*}
Since $\gamma \le \frac{C_2}{C_2+C_1^{-1}}$ we have that
$\log(\gamma^{-1})\geq (2C_1C_2)^{-1}$ by Maclaurin expansion 
and therefore
\begin{equation*}
\gamma^k=\exp(-k\log(\gamma^{-1}))\leq \exp(-k (2C_1C_2)^{-1}).
\end{equation*}
\end{proof}

In the last technical lemma we show that the error $\QQ \vv-\QQ^k\vv$ decays exponentially in $k$.
\begin{lemma}\label{lem:decay2}
For  any $\vv\in \VV$ it holds
\begin{equation*}
  |(\QQ - \QQ^k)\vv|_{\KK} \leq C_{\alpha,\beta,d,\mu,\sigma}k^{d/2}\exp(-k (6C_1C_2)^{-1})|\vv|_{\KK}.
\end{equation*}
\end{lemma}
\begin{proof}
  We use $\Rdd_T^\ell$ as an intermediate step to show that
  \begin{equation}\label{eq:decay}
    \begin{aligned}
      |\QQ_T \vv|_{\LL,\Omega \setminus U_{3\ell+1}(T)} & \le |\QQ_T \vv - \Rdd_T^\ell \vv|_{\LL} + | \Rdd_T^\ell \vv |_{\LL,\Omega \setminus U_{3\ell+1}(T)} \\
      &\le \alpha^{-1/2}|\QQ_T \vv - \Rdd_T^\ell \vv|_{\KK} + | \Rdd_T^\ell \vv |_{\LL,\Omega \setminus U_{3\ell+1}(T)} \\
&\le \alpha^{-1/2}\exp(-\ell (2C_1C_2)^{-1})|\vv|_{\KK,T}.
    \end{aligned}
\end{equation}
The term $ | \Rdd_T^\ell \vv |_{\LL,\Omega \setminus U_{3\ell+1}(T)}$ is zero since by equation (\ref{eq:spread}) $\Rdd_T^\ell \vv$ is zero outside $U_{3\ell}(T)$ and $\LL$ only spreads information one layer.

We let $\eta\in V_H$ be a cut-off function such that $(1-\eta)(x_i)=0$ for all $x_i\in U_{k-2}(T)$ and $(1-\III)(\eta\QQ_T \vv)\in \WW(U_k(T))$. Since $\QQ_T^k \vv$ is best approximation of $\QQ_T \vv$ in $\WW(U_k(T))$ we get 
\begin{align*}
|\QQ_T \vv-\QQ_T^k \vv|_\KK^2&\leq |\QQ_T \vv-(1-\III)(\eta \QQ_T \vv)|_\KK^2=|(1-\III)(\QQ_T \vv-\eta \QQ_T \vv)|_\KK^2
\\&\leq \beta |(1-\III)(\QQ_T \vv-\eta \QQ_T \vv)|_\LL^2 \leq C_{\beta,d,\mu,\sigma}|(1-\eta) \QQ_T \vv|_\LL^2
\end{align*}
using the equivalence of the $\LL$ and $\KK$ norms and Lemma~\ref{lem:interpolation}. Next we use the inequality \eqref{eq:product_rule}, since $\eta\in V_H$, 
and Lemma~\ref{lem:interpolation} to get
\begin{align*}
 |(1-\eta) \QQ_T \vv|_\LL^2&=\sum_{T'\in\mathcal{T}_H}  |(1-\eta) \QQ_T \vv|_{\LL,T'}^2\\
&=\sum_{T'\in\mathcal{T}_H}  |(1-\eta) \QQ_T \vv|_{\LL,T'\cap\Omega\setminus U_{k-3}(T)}^2\\
 &\leq 2\sum_{T'\in\mathcal{T}_H}  H^{-2}|\QQ_T \vv|_{\MM,T'\cap \Omega\setminus U_{k-3}(T)}^2+|\QQ_T \vv|_{\LL,T'\cap \Omega\setminus U_{k-3}(T)}^2\\
&\leq 2\sum_{T'\in\mathcal{T}_H}  H^{-2}|(1-\III)\QQ_T \vv|_{\MM,T'\cap \Omega\setminus U_{k-3}(T)}^2+|\QQ_T \vv|_{\LL,T'\cap \Omega\setminus U_{k-3}(T)}^2\\
 &\leq C_{d,\mu,\sigma}\sum_{T'\in\mathcal{T}_H} |\QQ_T \vv|_{\LL,T'\cap \Omega\setminus U_{k-6}(T)}^2\\
 &=C_{d,\mu,\sigma}|\QQ_T \vv|_{\LL,\Omega\setminus U_{k-6}(T)}^2.
\end{align*}
We use equation (\ref{eq:decay}) with $\ell=k/3-7/3$ to conclude
\begin{equation}\label{eq:QQT}
|\QQ_T \vv-\QQ_T^k \vv|_\KK\leq C_{\alpha,\beta,d,\mu,\sigma}\exp(-k (6C_1C_2)^{-1})|\vv|_{\KK,T}.
\end{equation}
Next we follow the proof of Theorem~4.3 in \cite{bookLOD}. Let $\eta\in V_H$ be $1$ for  $x\in \Omega\setminus U_{k+3}(T)$ and $0$ for all $x\in U_{k+2}(T)$. We let $\mathbf{e}=(\QQ -\QQ^k) \vv$ and $\mathbf{e}_T=(1-\III)(\eta \mathbf{e})\in \WW$ with $\mathbf{e}_T(x)=0$ for all $x\in U_{k+1}(T)$. We note that
$$
(\KK \mathbf{e}_T,\mathbf{e})=(\KK \mathbf{e}_T,\QQ_T \vv)=(\KK_T \mathbf{e}_T,\vv)=0
$$
since $\QQ_T^k\vv(x)=0$ for all $x\in\Omega\setminus U_k(T)$ and $\KK\mathbf{e}_T(x)=0$ for all $x\in U_k(T)$ and that $\KK_T \mathbf{e}_T=0$. We have $\mathbf{e}-\mathbf{e}_T=(1-\eta)\mathbf{e}+\III(\eta \mathbf{e})=(1-\III)((1-\eta)\mathbf{e})\in \WW(U_{k+4}(T))$. We conclude, using equation (\ref{eq:QQT}),
\begin{align*}
|\mathbf{e}|^2_\KK&=\sum_{T\in\mathcal{T}_H}(\KK(1-\III)((1-\eta)\mathbf{e}),(\QQ_T-\QQ_T^k)\vv)\\
&\leq C_{\alpha,\beta,d,\mu,\sigma}\text{exp}(-k(6C_1C_2)^{-1})\sum_{T\in\mathcal{T}_H}|\mathbf{e}|_{\KK,U_{k+5}(T)}|\vv|_{\KK,T}\\
&\leq C_{\alpha,\beta,d,\mu,\sigma}k^{d/2}\text{exp}(-k(6C_1C_2)^{-1})|\mathbf{e}|_\KK |\vv|_{\KK}
\end{align*}
where we use that one element $T$ is only covered by a finite number of patches $U_{k+5}(T)$. The lemma follows after division by $|\mathbf{e}|_\KK$.
\end{proof}

With Lemma~\ref{lem:decay2}, the use of $\QQ^k$ instead of $\QQ$ is thoroughly motivated. Moreover, with exponential decay the element patches can be small and still be representative. Now all that is left is to provide the final a priori estimate for the localized LOD approximation $\uu^k_H$:

\begin{theorem}\label{theorem:bound} 
  Under Assumptions~\ref{ass:K} and \ref{ass:network} with $H\geq 4dR_0$ the error in the approximate solution $\uu^k_H$, defined in equation (\ref{eq:u_Hk}), fulfills
  $$
  |\uu-\uu_H^{k}|_\KK\leq C_{\alpha,\beta,d,\mu,\sigma}k^{d/2}\left(H +\exp(-k (6C_1C_2)^{-1})\right)|\ff|_{\MM^{-1}}.
  $$
  \end{theorem}
  \begin{proof}
  By Galerkin orthogonality $|\uu-\uu_H^k|_\KK\leq |\uu-\vv|_\KK$ for all $\vv\in \VV_H^{\text{ms},k}$. We let $\vv=(1-\QQ^k)\III\uu_H\in \VV_H^{\text{ms},k}$ and use the identity $\uu_H=(1-\QQ)\III\uu_H$. Using the triangle inequality we therefore have
  $$
  |\uu-\uu_H^k|_\KK\leq |\uu-\uu_H|_\KK+|(\QQ-\QQ^k)\III\uu_H|_\KK.
  $$
  The first part is treated in Lemma~\ref{lem:ideal}. For the second part we use the triangle inequality, Lemma~\ref{lem:decay2}, and that $\LL$ and therefore $\KK$ is stable with respect to $\III$ in the $\KK$ norm
  \begin{align*}
  |(\QQ-\QQ^k)\III \uu_H|_\KK^2  &\leq C_{\alpha,\beta,d,\mu,\sigma}k^{d}\exp(-k (3C_1C_2)^{-1})|\III \uu_H|_\KK^2\\
  &\leq C_{\alpha,\beta,d,\mu,\sigma}k^{d}\exp(-k (3C_1C_2)^{-1})|{\uu_H}|_\KK^2\\
  &\leq C_{\alpha,\beta,d,\mu,\sigma}k^{d}\exp(-k (3C_1C_2)^{-1})|\ff|_{\MM^{-1}}^2.
  \end{align*}
  where we use the equivalence of the $\KK$ and $\LL$ norms and that $|\uu_H|_\KK\leq |\uu|_\KK\leq C_\alpha|\ff|_{\MM^{-1}}$. The theorem follows.
\end{proof}

\section{Numerical examples}\label{sec:num}
\begin{figure}
  \centering
  \includegraphics[width = 0.25\textwidth]{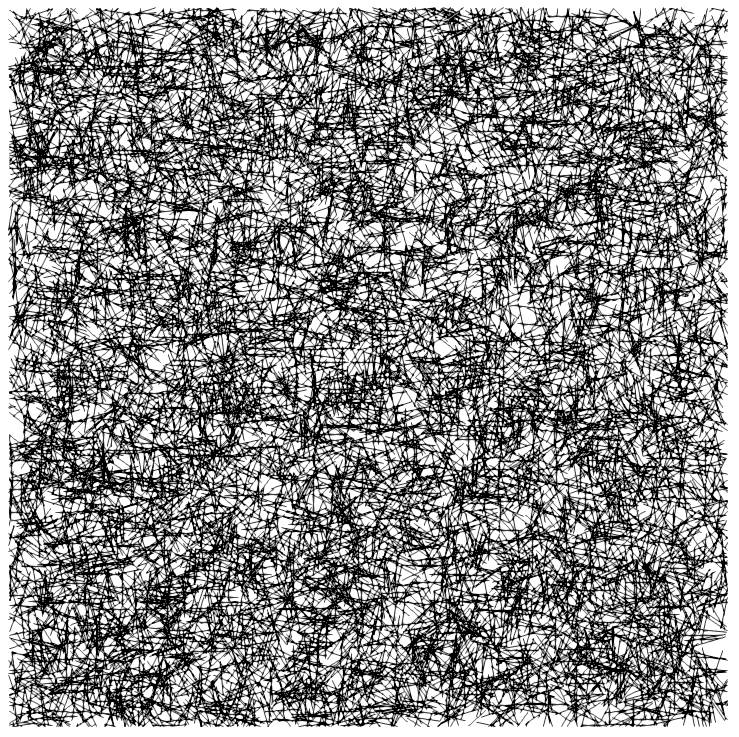}
  \caption{The network analyzed in the numerical examples. It is constructed by around $20000$ line segments of length $0.05$ uniformly distributed  in the unit square. All nodes sits at crossings of line segments. Edges are placed between two nodes that share a common line segment. }
  \label{fig:network}
\end{figure}

\noindent 
We first consider a scalar example modelling heat conduction and then we turn to a structural problem where we seek the displacement of a fiber network. For all numerical examples, we use the network shown in Figure~\ref{fig:network}. Uniformly rotated line segments of length $0.05$ are uniformly distributed in the unit square so that the total mass is $|1|_M^2=1000$, resulting in about $20000$ line segments. The network nodes are then defined as the line segments' endpoints and crossings, with the network edges connecting every two nodes that share a line segment. The total number of nodes in the generated network is around 80000.  
\subsection{Heat conduction} \label{num:heat}
\noindent We consider the model problem introduced in  Example  \ref{example:heat} for the two-dimensional network in Figure~\ref{fig:network} and adopt the scalar notation from Example  \ref{example:heat}. The solution represents temperature (scalar) in each node, and the node-wise operator, $K_x$, is defined as:
$$(K_x v,v) = \frac{1}{2} \sum_{y\sim x}\gamma_{xy} \frac{(v(x)-v(y))^2}{|x - y|},$$
where the coefficients $\gamma_{xy}\in[0.1,1]$ are chosen at random for each edge $\{x,y\}$. The computational domain is the unit square $\Omega=[0,1]^2$. The problem considered has a constant right hand side weighted with the mass matrix $M$ and zero Dirichlet boundary is applied on the entire boundary, i.e.
\begin{equation*}
  \begin{cases}
  K u = f, \\
  u(\partial \Omega)  = 0,
  \end{cases}
\end{equation*}
with $K=\sum_{x\in\mathcal{N}} K_x$, $f=M1$ and $1\in \hat{V}$. The exact solution is compared to the LOD approximation \eqref{eq:loclod} with localization parameter $k=2$ for different coarse grids. To show that the problem can not easily be solved using the coarse finite element spaces $V_H$ we also consider the problem:
\begin{equation}
  \text{find } u^{\text{FEM}}_H\in V_H \,:\, (K u^{\text{FEM}}_H, v)  = (M 1, v) \text{ for all } v \in V_H.
\end{equation}

\begin{figure}
  \centering
  \includegraphics[width = 0.3\textwidth]{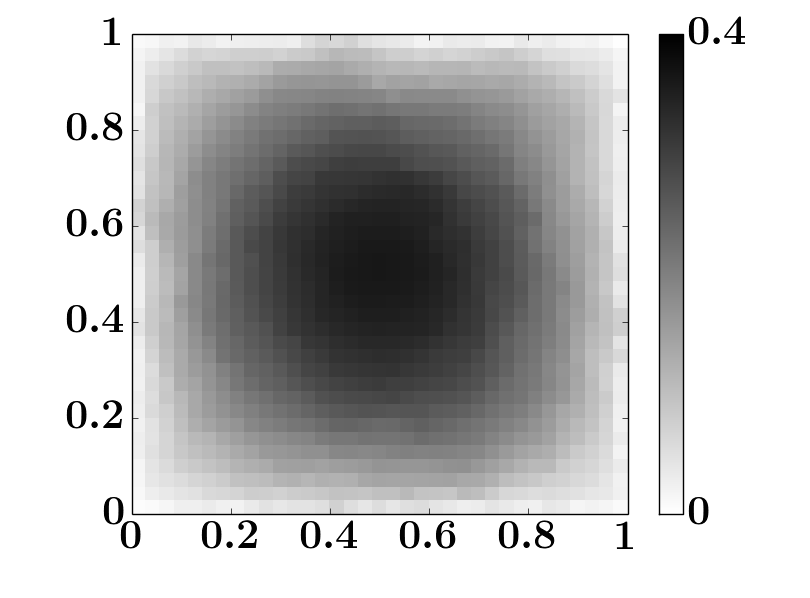}
  \includegraphics[width = 0.3\textwidth]{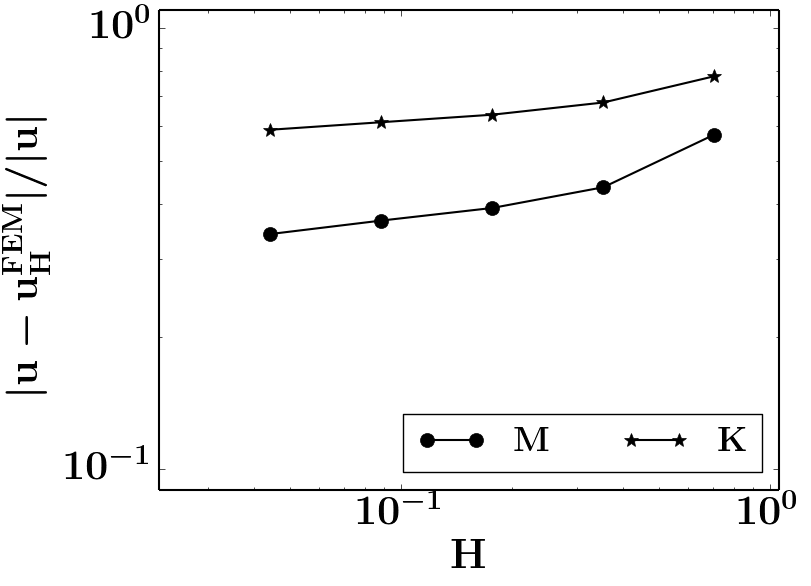}
  \includegraphics[width = 0.3\textwidth]{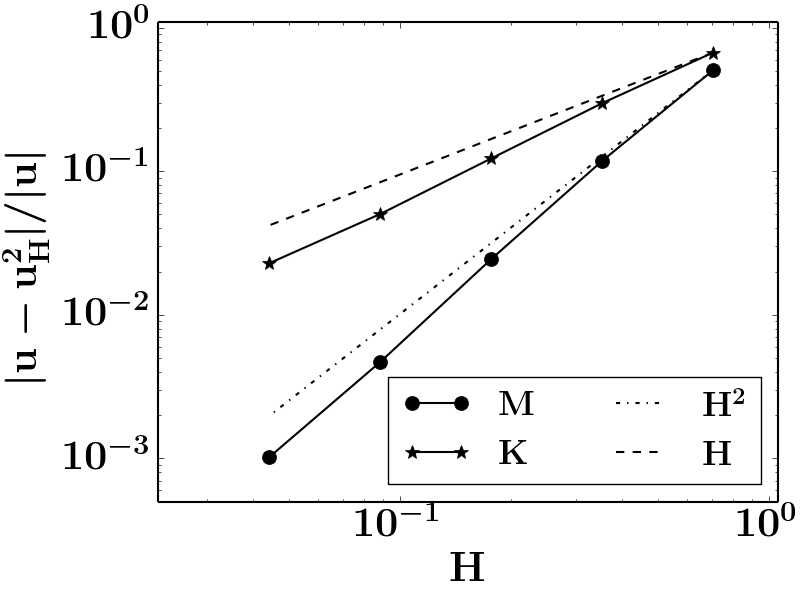}
  \caption{The solution $u$ to the problem in section \ref{num:heat}, along with the convergence results for a finite element approach and the LOD method with localization factor $k=2$.}
  \label{fig:heatnum}
\end{figure}

An illustration of the solution, $u$, and the errors of the direct finite element approach and the LOD approximations in both $K$ and $M$ norm are presented in Figure~\ref{fig:heatnum}. The results show a convergence plateau for the finite element approach, whereas the theoretical convergence rate of $H$ (Theorem~\ref{theorem:bound}) is achieved for the LOD method already for a localization parameter of $k=2$. Moreover, we observe that error is proportional to $H^2$ in the $M$-norm for the LOD method.

\subsection{A fiber network model}\label{num:fiber}
Here two variations of Example~\ref{ex:fiber} are considered. The network in Figure~\ref{fig:network} should be interpreted as a mesh of round steel wires of radius $r_w = 2.5\times 10^{-3}$. Equation \eqref{eq:spring} is a linearized version of Hooke's law with parameter $\gamma_{xy} =\gamma_1 = EA$, where $A = \pi r_w^2$ is the cross-section area of the wire and $E = 210$\,GPa its Young's modulus. The bending forces are handled by adding the equations in \eqref{eq:fiber}. These additions are linearized versions of Euler--Bernoulli with parameters $\gamma_{xyz} = 2EI(|x-y|+|x-z|)^{-2}$ where $E$ is the same Young's modulus and $I = 0.25\pi r_w^4 = 0.25Ar_w^2$ is the second moment of area of the wire. The two coefficients are related in the following way,
$$\gamma_{xyz} = EA\frac{r_w^2}{2(|x-y|+|x-z|)^{2}} = \gamma_{1}\frac{r_w^2}{2(|x-y|+|x-z|)^{2}},$$
where $0.05 \leq \frac{r_w}{|x-y|}\leq 5$ for any edge $x\sim y$. This relation is dependent on the lengths of the edges, where the edge lengths in turn depend on how the fibers intersect each other. Because of this $\gamma_{xyz}$ varies rapidly in space.

\subsubsection*{Pure displacement problem}
\begin{figure}
  \centering
  \hspace{0.1\textwidth}
  \includegraphics[width = 0.40\textwidth]{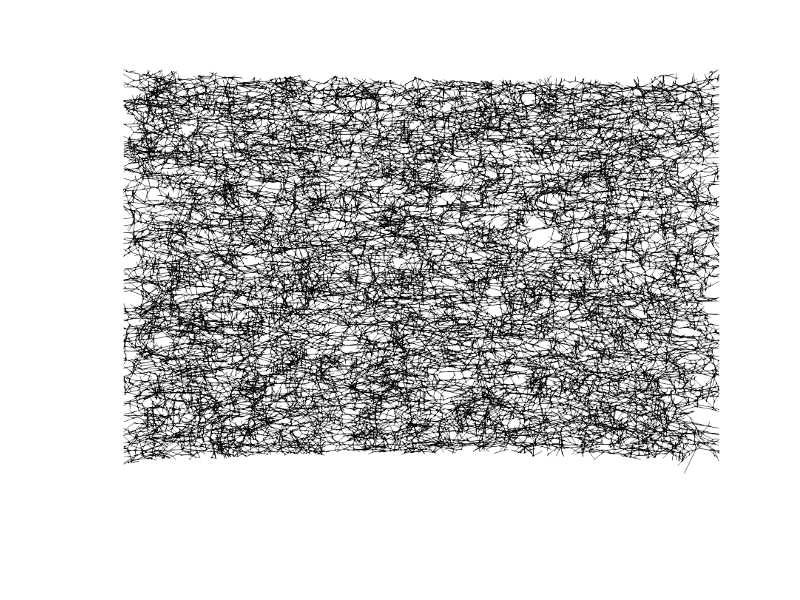}
  \includegraphics[width = 0.40\textwidth]{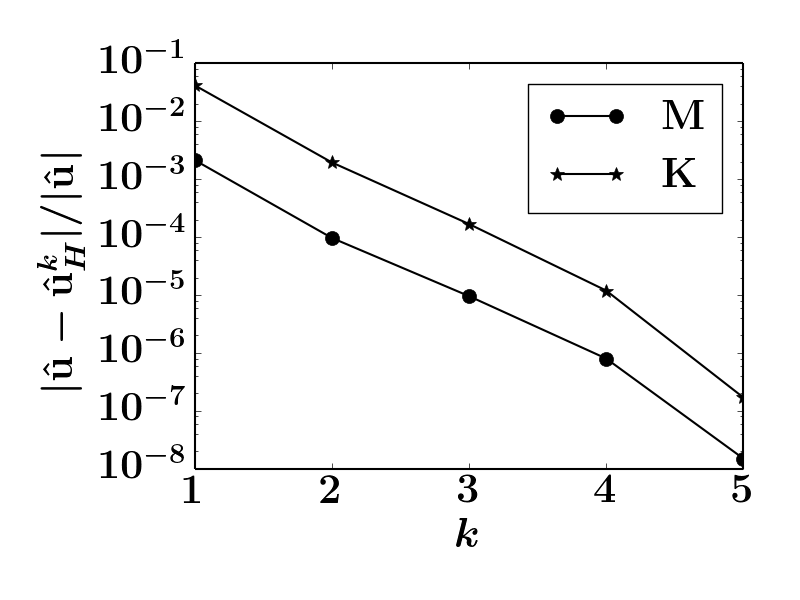}
  \caption{The solution of the strained fiber network along with the normalized approximation errors for different localization paramters $k$.}
  \label{fig:disp}
\end{figure}

The first structural problem we consider is a tensile simulation, where one side of the unit square is fixed, and the opposite side is displaced. This displacement stretches the network, and the solution to the problem is the equilibrium of the network given the displacement. We will only consider forces and displacements in the plane the network resides in for this simulation, meaning that any $x_3$-directional components are left out. The problem can be written as 
\begin{equation*}
  \begin{cases}
  \KK \hat{\uu}  = \mathbf{0}, \\
  \hat{\uu} (\Gamma_1)  = [0,0]^T, \ \hat{\uu} (\Gamma_2)  = [0.5,0]^T, \\
  \end{cases}
\end{equation*}
where $\Gamma_1$ is any point with x-coordinate 0, and $\Gamma_2$ is any point with $x$-coordinate 1. The solution is presented in Figure~\ref{fig:disp}. Solving this problem with the LOD method requires some extra steps compared to the previous example as we have non-vanishing Dirichlet conditions. As mentioned in Section \ref{section:problem}, we introduce an auxiliary function, $\mathbf{g}$, such that $\hat{\uu}=\uu+\mathbf{g}$ and consider,
\begin{equation*}
  \begin{cases}
    \KK \uu =  -\KK \mathbf{g}, \\
    \uu (\Gamma_1\cup\Gamma_2)  = [0,0]^T.
  \end{cases}
\end{equation*}
For this specific problem we choose $\mathbf{g}(x) = [0.5x_1,0]^T$ which is in $\hat{\mathbf{V}}_H$ for all $H$. It was shown in \cite{BIT} that if $\mathbf{g}\in\hat{\mathbf{V}}_H$  then the exact solution of this support problem can be written as $\uu=\uu_H + \mathbf{c}_H$, where $\mathbf{c}_H$ is attainable with an extended version of $\QQ$. This is seen by first writing the corrector term, $\mathbf{c}_H$, as the solution to the following variational problem:
$$\text{find } \mathbf{c}_H\in \WW \,:\, (\KK\mathbf{c}_H, \ww) = (\KK (-\mathbf{g}), \ww)\text{ for all } \ww \in \WW,$$
by using that $\mathbf{V} = \mathbf{V}^{\text{ms}} \bigoplus \mathbf{W}$, $\mathbf{V}^{\text{ms}}\perp_{\mathbf{K}} \mathbf{W}$, and $\KK$ being coercive. The solution to this variational problem can be written as $\hat{\QQ}(-\mathbf{g})$, where $\hat{\QQ} : \hat \VV \rightarrow \WW$ is the extended projection operator of $\mathbf{Q}$:
$$(\KK \hat{\mathbf{Q}} \mathbf{v},\ww) =  (\KK \mathbf{v},\ww) \text{ for all } \mathbf{w}\in \mathbf{W}.$$
With $\hat{\QQ}_T$ and $\hat{\QQ}_T^k$ derived analogously to $\QQ_T$ and $\QQ_T^k$. In practice, finding this extension, $\hat \QQ$, is comparable to finding $\QQ$ in terms of computational complexity. Using this projection operator we can write the exact solution to the initial problem as: 
$$\hat{\uu} = \uu+\mathbf{g} = \uu_H+ \mathbf{c}_H + \mathbf{g} =  \uu_H+(1-\hat{\QQ})\mathbf{g},$$
and the localized LOD approximations:
$$\hat{\uu}_H^k = \uu_H^k+(1-\hat{\mathbf{Q}}^k)\mathbf{g}.$$
\noindent For an extended discussion on how to handle general Dirichlet data in the LOD method, see \cite{generalDirichlet}.

In the numerical experiment the exponential decay of the correctors are analyzed, by fixing $H=1/32$ and computing the errors $|\hat{\uu}^k_H-\hat{\uu}|$ for different $k$. The results are presented in Figure~\ref{fig:disp}, and exponential decay is observed in both the $\KK$ and $\MM$ norm which is consistent with Theorem~\ref{theorem:bound}.

\subsubsection*{Displacement problem with lateral load}
\begin{figure}
  \centering
  \hspace{0.1\textwidth}
  \includegraphics[width = 0.40\textwidth]{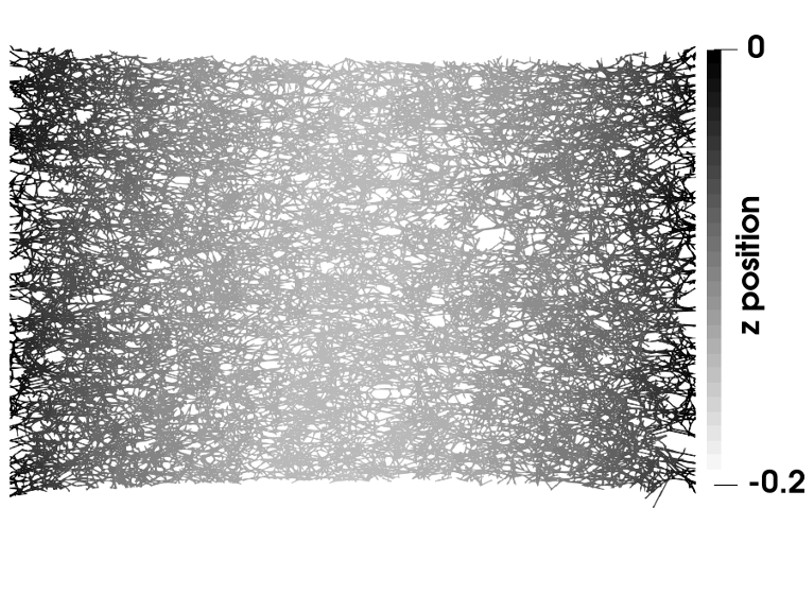}
  \includegraphics[width = 0.40\textwidth]{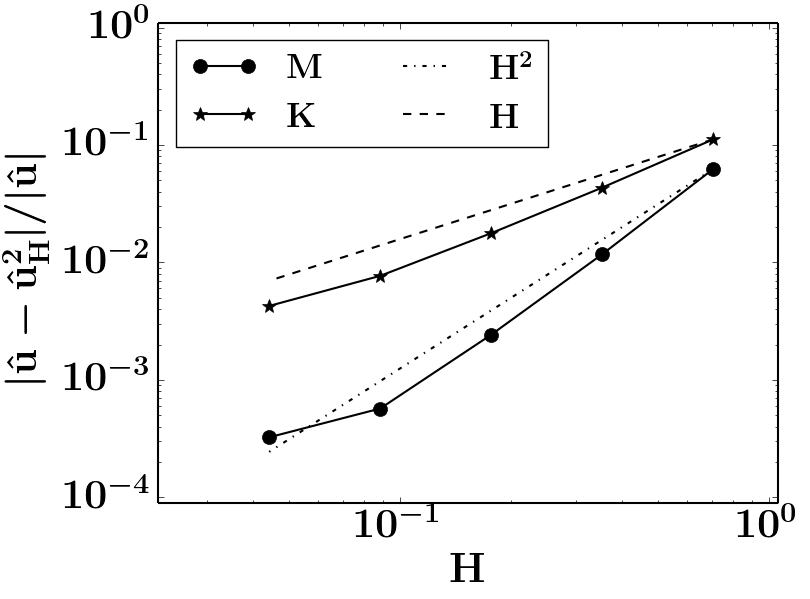}
  \caption{Illustration of the solution (XY-plane) and the normalized errors of the LOD approximations $\hat{\uu}^2_H$ for varying $H$.}
  \label{fig:disp3d}
\end{figure}

In the second numerical example of the fiber network problem, we introduce a lateral ($x_3$-directional) load to the previous tensile problem. This problem can be expressed as the following linear system:
\begin{equation*}
  \begin{cases}
  \KK \hat{\uu} = \ff, \\
  \hat{\uu} (\Gamma_1)  = [0,0,0]^T, \ \hat{\uu} (\Gamma_2)  = [0.5,0,0]^T, \\
  \end{cases}
\end{equation*}
where $\Gamma_1$ is any point with x-coordinate $0$, $\Gamma_2$ any point with x-coordinate $1$, and $\ff=\MM \mathbf{h}$ with $\mathbf{h}$ as the constant function $[0,0,-10^5]$. As with the previous example, we let $\mathbf{g}(x) = [0.5x_1,0,0]^T \in \hat{\VV}_H$ and 
\begin{equation*}
  \begin{cases}
  \KK \uu = \ff-\KK \mathbf{g}, \\
  \uu (\Gamma_1\cup\Gamma_2)  = [0,0,0]^T
  \end{cases}
\end{equation*}
where $\hat{\uu} = \uu+\mathbf{g}$. Using the same motivation as in the previous example, the localized LOD approximations considered are:
$$\hat{\uu}_H^k =  \uu_H^k+(1-\hat{\mathbf{Q}}^k)\mathbf{g}.$$ 
However, unlike the previous example we can not guarantee that the ideal LOD approximation, $\hat{\uu}_H$, is the exact solution $\hat{\uu}$, since $\ff\neq \mathbf{0}$. 
Theorem~\ref{theorem:bound} is numerically confirmed, with localization parameter $k=2$, and presented in Figure~\ref{fig:disp3d}, along with the reference solution $\hat \uu$. The $H$ convergence is seen in the $\KK$-norm as the theory indicates, but some slight stagnation is observed for the finest grid considered which would vanish for $k=3$ as indicated in Figure~\ref{fig:disp}. 
Already for $k=2$ the method produces highly accurate results in both $\KK$ and $\MM$ norm, with less than one percent relative error in the $\KK$ norm and a tenth of a percent in the $\MM$ norm for the finest coarse grid considered. 

\section*{Acknowledgements}
The Swedish Foundation for Strategic Research (SSF) supports the second author, the third and fifth author are supported by the G\"oran Gustafsson Foundation for Research in Natural Sciences and Medicine, and the Swedish Research Council supports the fourth and fifth author in project number 2019-03517 VR.

\bibliographystyle{abbrv}

\begin{thebibliography}{99}

  
\bibitem{multigrid} Brandt, A., {\em Multi-level adaptive solutions to boundary-value problems,} Mathematics of Computation, 31, pp. 333–390, 1977.


\bibitem{Chu} Chu, J., Engquist, B., Prodanovi\'{c}, M., and Tsai, R., {\em A multiscale method coupling network and continuum models in porous media I: steady-state single phase flow,} Multiscale Model. Simul. 10, pp. 515-549, 2012.

\bibitem{Rossa} Della Rossa, F., D’Angelo, C., and Quarteroni, F, {\em A distributed model of traffic flows on extended regions,} Netw. Heterog. Media 5, pp. 525-544, 2010.

\bibitem{msfem} Efendiev, Y., Galvis, J., and Hou, T.Y., {\em Generalized multiscale finite element methods (GMsFEM),} Journal of Computational Physics archive, 251, pp. 116-135, 2013.


\bibitem{Ewing} Ewing, R., Iliev, O., Lazarov, R., Rybak, I., and Willems, J., {\em A simplified method for upscaling composite materials with high contrast of the conductivity,} SIAM J. Sci. Comput., 31, pp. 2568-2586, 2009.

  
  

\bibitem{GoHeMa22} G\"{o}rtz, M., Hellman, F., and M\aa lqvist, A., {\em Iterative solution of spatial network models by subspace decomposition,} preprint arXiv:2207.07488 

\bibitem{BIT} Kettil, G., M\aa lqvist, A., Mark, A., Fredlund, M., Wester, K., and Edelvik, F., {\em Numerical upscaling of discrete network models,} BIT, 60, pp. 67-92, 2020.

\bibitem{generalDirichlet}  Henning, P.~and Målqvist, A., { \em Localized orthogonal decomposition techniques for boundary value problems}, SIAM J. Sci. Comp., 36, A1609-A1634, 2014.

\bibitem{Iliev} Iliev, O., Lazarov, R., and Willems, J., {\em Fast numerical upscaling of heat equation for fibrous materials,} Comput. Visual. Sci., 13, pp. 275-285, 2010.

\bibitem{KPY18} Kornhuber, R., Peterseim, D., and Yserentant, H., {\em An analysis of a class of variational multiscale methods based on subspace decomposition,} Mathematics of Computation, 87(314):2765-2774, 2018.

\bibitem{KY16} Kornhuber, R.~and Yserentant, H., {\em Numerical homogenization of elliptic multiscale problems by subspace decomposition,} Multiscale Model. Simul., 14(3), pp. 1017-1036, 2016. 

\bibitem{aprioriLOD} M\aa lqvist, A.~and Peterseim, D., {\em Localization of elliptic multiscale problems,} Math. Comp. 83, pp. 2583-2603, 2014.

\bibitem{bookLOD} M\aa lqvist, A.~and Peterseim, D., {\em Numerical homogenization by localized orthogonal decomposition,} SIAM Spotlights, ISBN: 978-1-611976-44-1, 2020  

\bibitem{gamblets} Owhadi, H.~and Scovel, C., {\em Operator-adapted wavelets, fast solvers, and numerical homogenization,} volume 35 of {\em Cambridge Monographs on Applied and Computational Mathematics,} Cambridge University Press, Cambridge, UK, 2019. 

\bibitem{SZ90} Scott, R.~and Zhang, S., {\em Finite element interpolation of nonsmooth functions satisfying boundary conditions,} Math. Comp., 54, pp. 483-493, 1990.

\bibitem{paper}
Svenning, E., Mark, A., Edelvik, F., Glatt, E., Rief, S., Wiegmann, A., Martinsson, L., Lai, R., Fredlund, M.~and Nyman, U, {\em Multiphase simulation of fiber suspension flows using immersed boundary methods,} Nordic Pulp and Paper Research Journal, 27(2), pp. 184--191, 2012.

\end{thebibliography}

\end{document}